\documentclass[leqno,12pt]{article}

\usepackage{amssymb}
\usepackage{euscript}
\usepackage{natbib}

\evensidemargin\oddsidemargin

\begin{document}

\baselineskip=18pt \setcounter{page}{1}

\renewcommand{\theequation}{\thesection.\arabic{equation}}
\newtheorem{theorem}{Theorem}[section]
\newtheorem{lemma}[theorem]{Lemma}
\newtheorem{proposition}[theorem]{Proposition}
\newtheorem{corollary}[theorem]{Corollary}
\newtheorem{remark}[theorem]{Remark}
\newtheorem{fact}[theorem]{Fact}
\newtheorem{problem}[theorem]{Problem}

\newcommand{\eqnsection}{
\renewcommand{\theequation}{\thesection.\arabic{equation}}
    \makeatletter
    \csname  @addtoreset\endcsname{equation}{section}
    \makeatother}
\eqnsection

\def\r{{\mathbb R}}
\def\e{{\mathbb E}}
\def\p{{\mathbb P}}
\def\bZ{{\mathbb Z}}
\def\S{{\mathbb S}}
\def\t{{\mathbb T}}
\def\bN{{\mathbb N}}
\def\deg{{\rm b}}
\def\P{{\bf P}}
\def\ind {\mbox{\rm 1\hspace {-0.04 in}I}}
\def\E{{\bf E}}
\def\ee{\mathrm{e}}
\def\d{\, \mathrm{d}}
\newcommand{\cF}{\mathcal{F}}
\newcommand{\cG}{\mathcal{G}}
\newcommand{\cU}{\mathcal{U}}
\def\Q{{\bf Q}}
\def\q{{\mathbb Q}}


\vglue15pt

\centerline{\Large\bf Transient random walks in random environment}

\bigskip

\centerline{\Large\bf on a Galton--Watson tree}
\bigskip
\bigskip
\bigskip
\centerline{by}

\bigskip

\centerline{Elie Aid\'ekon}
\medskip

\centerline{\it Universit\'e Paris VI}

\bigskip
\bigskip

\medskip

\bigskip
\bigskip
\bigskip
\bigskip

{\leftskip=2truecm \rightskip=2truecm \baselineskip=15pt \small

\noindent{\slshape\bfseries Summary.} We consider a transient random
walk $(X_n)$ in random environment on a Galton--Watson tree. Under
fairly general assumptions, we give a sharp and explicit criterion
for the asymptotic speed to be positive. As a consequence,
situations with zero speed are revealed to occur. In such cases, we
prove that $X_n$ is of order of magnitude $n^{\Lambda}$, with
$\Lambda \in (0,1)$. We also show that the linearly edge reinforced
random walk on a regular tree always has a positive asymptotic
speed, which improves a recent result of Collevecchio~\cite{Col06}.

\bigskip

\bigskip


\bigskip

} 

\bigskip
\bigskip

\section{Introduction}
   \label{s:intro}

\subsection{Random walk in random environment}

Let $\nu$ be an ${\bN^*}$-valued random variable (with ${\bN^*} :=
\{ 1,2, \cdots\}$) and $(A_i, i\ge 1)$ be a random variable taking
values in $\r_+^{\bN^*}$. Let $q_k:=P(\nu=k)$, $k\in {\bN^*}$. We
assume $q_0=0$, $q_1<1$, and $m:=\sum_{k\ge 0}kq_k<\infty$. Writing
$V:= (A_i,i\le \nu)$, we construct a Galton--Watson tree as follows.

Let $e$ be a point called the root. We pick a random variable
$V(e):=(A(e_i),i\le \nu(e))$ distributed as $V$, and draw $\nu(e)$
children to $e$. To each child $e_i$ of $e$, we attach the random
variable $A(e_i)$. Suppose that we are at the $n$-th generation. For
each vertex $x$ of the $n$-th generation, we pick independently a
random vector $V(x)=(A(x_i),i\le \nu(x))$ distributed as $V$,
associate $\nu(x)$ children $(x_i,i\le \nu(x))$ to $x$, and attach
the random variable $A(x_i)$ to the child $x_i$. This leads to a
Galton--Watson tree $\t$ of offspring distribution $q$, on which
each vertex $x\neq e$ is marked with a random variable $A(x)$.

We denote by $GW$ the distribution of $\t$. For any vertex $x\in
\t$, let $\buildrel \leftarrow  \over x$ be the parent of $x$ and
$|x|$ its generation ($|e|=0$). In order to make the presentation
easier, we artificially add a parent $\buildrel \leftarrow  \over e$
to the root $e$. We define the environment $\omega$ by
$\omega(\buildrel \leftarrow  \over e,e)=1$ and for any vertex $x\in
\t \backslash \{ {\buildrel \leftarrow \over e} \} $,
\begin{itemize}
\item $\omega(x,x_i)={A(x_i)\over 1+\sum_{i=1}^{\nu(x)}A(x_i)}$\,, $\forall\, 1\le i\le \nu(x)$,
\item $\omega(x,\buildrel \leftarrow  \over x)={1\over 1+\sum_{i=1}^b A(x_i)}$\,.
\end{itemize}

\noindent For any vertex $y\in\t$, we define on $\t$ the Markov
chain $(X_n, \, n\ge 0)$ starting from $y$ by
\begin{eqnarray*}
P_{\omega}^y(X_0=y) &=&1,\\
P_{\omega}^y(X_{n+1}=z \, | \, X_n=x)&=&\omega(x,z)\,.
\end{eqnarray*}

\noindent Given $\t$, $(X_n, \, n\ge 0)$ is a $\t$-valued random
walk in random environment (RWRE). We note from the construction
that $\omega(x,.)$, $x\neq \buildrel \leftarrow \over e$ are
independent.

Following \cite{lp92}, we also suppose that $A(x)$, $x\in\t$,
$|x|\ge 1$, are identically distributed. Let $A$ denote a random
variable having the common distribution. We assume the existence of
$\alpha
>0$ such that $\hbox{ess sup} (A)\le \alpha$ and $\hbox{ess sup} ({1\over A})\le \alpha$. The following criterion is known. \\

\noindent {\bf Theorem A (Lyons and Pemantle~\cite{lp92})} {\it The walk $(X_n)$ is transient if $\,\inf_{[0,1]} \E[A^t]> \frac{1}{m}$, and is recurrent otherwise.}\\

When $\t$ is a regular tree, Menshikov and Petritis~\cite{mp02}
obtain the transience/recurrence criterion by means of a
relationship between the RWRE and Mandelbrot's multiplicative
cascades; Hu and Shi~\cite{hs06},\cite{hs05} characterize  different
asymptotics of the walk in the recurrent case, revealing a wide
range of regimes.

Throughout the paper, we assume that the walk is transient (i.e.,
$\inf_{[0,1]} \E[A^t]> \frac{1}{m}$ according to Theorem A). Given
the transience, natural questions arise concerning the rate of
escape of the walk. The law of large numbers says that there exists
a deterministic $v\ge 0$ (which can be zero) such that
$$
\lim_{n\rightarrow \infty} \frac{|X_n|}{n}=v, \qquad a.s.
$$

\noindent This was proved by Gross~\cite{gross} when $\t$ is a
regular tree, and by Lyons et al.~\cite{lpp96} when $A$ is
deterministic; their arguments can be easily extended in the general
case (i.e., when $\t$ is a Galton--Watson tree and $A$ is random).

We are interested in determining whether $v>0$.

When $A$ is deterministic, it is shown by Lyons et al.~\cite{lpp96}
that the transient random walk always has positive speed. Later, an
interesting large deviation principle is obtained in Dembo et
al.~\cite{dgpz02}. In the special case of non-biased random walk,
Lyons et al.~\cite{lpp95} succeed in computing the value of the
speed.

We recall two results for RWRE on $\bZ$ (which can be seen as a half
line-tree). The first one gives a necessary and sufficient condition
for RWRE to have positive asymptotic speed.

\bigskip

\noindent {\bf Theorem B (Solomon~\cite{solomon75})} {\it If
$\,\t=\bZ$, then
$$
\E\left[{1\over A}\right]<1 \Longleftrightarrow
\lim_{n\rightarrow\infty}{X_n\over n}>0 ~~ a.s.
$$}

When the transient RWRE has zero speed, Kesten, Kozlov and Spitzer
in \cite{Kes75} prove that the walk is of polynomial order. To this
end, let $\kappa\in(0,1]$ be such that $E\left[{1\over
A^{\kappa}}\right]=1$. Under some mild conditions on $A$,
\begin{itemize}
\item if $\kappa<1$, then ${X_n \over n^{\kappa}}$ converges in
distribution.
\item If $\kappa=1$, then ${\ln(n) X_n \over n}$
converges in probability to a positive constant.
\end{itemize}

\bigskip
\bigskip

The aim of this paper is to study the behaviour of the transient
random walk  when $\t$ is a Galton--Watson tree. Let $Leb$ represent
the Lebesgue measure on $\r$ and let
\begin{equation}
    \Lambda:=Leb \left\{t\in \r: \E[A^t]\le
    {1\over q_1}\right\} .
    \label{Lambda}
\end{equation}

\noindent If $q_1=0$, then we define $\Lambda:=\infty$. Notice that
this definition is similar to the definition of $\kappa$ in the
one-dimensional setting. Our first result, which is a (slightly
weaker) analogue of Solomon's criterion for Galton--Watson tree
$\t$, is stated as follows.

\begin{theorem}
\label{global}
 Assume $\inf_{[0,1]} \E[A^t]> \frac{1}{m}$, and let
 $\Lambda$ be as in $(\ref{Lambda})$.

 {\rm (a)} If $\Lambda<1$, the walk has zero speed.

 {\rm (b)} If $\Lambda>1$, the walk has positive
 speed.
\end{theorem}

\begin{corollary}
\label{c:global}
 Assume $\inf_{[0,1]} \E[A^t]> \frac{1}{m}$. If $\,\t$ is a regular tree, then the walk has positive speed.

\end{corollary}
\bigskip

Theorem \ref{global} extends Theorem B, except for the ``critical
case" $\Lambda=1$.

Corollary \ref{c:global} says there is no
Kesten--Kozlov--Spitzer-type regime for RWRE when the tree is
regular. Our next result exhibits such a regime for Galton--Watson
trees $\t$.

\begin{theorem}
\label{kesten}
 Assume $\inf_{[0,1]} \E[A^t]> \frac{1}{m}$, and  $\Lambda\le 1$. Then
$$
\lim_{n\rightarrow\infty}{\ln(|X_n|)\over \ln(n)}=\Lambda~~a.s.
$$
\end{theorem}
\bigskip

Since $\Lambda>0$, the walk is proved to be of polynomial order. As
expected, $\Lambda$ plays the same role as $\kappa$.

\subsection{Linearly edge reinforced random walk}

The reinforced random walk is a model of random walk introduced by
Coppersmith and Diaconis \cite{CoDi} where the particle tends to
jump to familiar vertices. We consider the case where the graph is a
$b$-ary tree $\t$, that is a tree where each vertex has $b$ children
($b\ge 2$). At each edge $(x,y)$, we initially assign the weight
$\pi(x,y)=1$. If we know the weights and the position of the walk at
time $n$, we choose an edge emanating from $X_n$ with probability
proportional to its weight. The weight of the edge crossed by the
walk then increases by a constant $\delta>0$. This process is called
the Linearly Edge Reinforced Random Walk (LERRW). Pemantle in
\cite{Pe88} proves that there exists a real $\delta_0$ such that the
LERRW  is transient if $\delta<\delta_0$ and recurrent if
$\delta>\delta_0$ ($\delta_0=4,29..$ for the binary tree). We focus,
from now on, on the case $\delta=1$, so that the LERRW almost surely
is transient. Recently, Collevecchio in \cite{Col06} shows that when
$b\ge 70$ the LERRW has a positive speed $v$ which verifies $0<v\le
{b\over b+2}$. We propose to extend the positivity of the speed to
any $b\ge 2$.

\begin{theorem}
\label{reinforce} The linearly edge reinforced random walk on a
$b$-ary tree has positive speed.
\end{theorem}

We rely on a correspondence between RWRE and LERRW, explained in
\cite{Pe88}. By means of a Polya's urn model, Pemantle shows that
the LERRW has the distribution of a certain RWRE, such that for any
$y\neq {\buildrel \leftarrow\over e}$, the density of $\omega(y,z)$
on $(0,1)$ is given by
\begin{itemize}
\item  $f_0(x)=\frac{b}{2}\,(1-x)^{\frac{b}{2}-1}$\qquad \qquad \qquad ~~ if $z={\buildrel \leftarrow \over y}$,
\item  $f_1(x)= {\Gamma(\frac{b}{2}\,+1)\over
\Gamma(\frac{1}{2})\Gamma(\frac{b+1}{2})}\,x^{\!-\frac{1}{2}}(1-x)^{\frac{b-1}{2}}$\qquad
if $z$ is a child of $y$.
\end{itemize}

\noindent Consequently, we only have to prove the positivity of the
speed of this RWRE.

With the notation of Section 1.1, $A$ is not bounded in this case,
which means Theorem \ref{global} does not apply. To overcome this
difficulty, we prove the following result.

\begin{theorem}
\label{errw} Let $\t$ be a $b$-ary tree and assume that
$\inf_{[0,1]} \E[A^t]> \frac{1}{b}$ and
$$E\left[\left(\sum_{i=1}^b A_i\right)^{-1}\right]<\infty\,.$$ Then
the RWRE has positive speed.
\end{theorem}

Since the RWRE associated with the LERRW satisfies the assumptions
of Theorem \ref{errw} as soon as $b\ge 3$, Theorem \ref{reinforce}
follows immediately in the case $b\ge 3$. The case of the binary
tree is dealt with separately.

\bigskip

The rest of the paper is organized as follows. We prove Theorem
\ref{errw} in Section 2. In Section 3, we prove the upper bound in
Theorem \ref{kesten}. Some technical results are presented in
Section 4, and are useful in Section 5 in the proof of the lower
bound in Theorem \ref{kesten}. In Section 6, we prove Theorem
\ref{global}. The proof of Theorem \ref{reinforce} for the binary
tree is the subject of Section 7. Finally, Section 8 is devoted to
the computation of parameters used in the proof of Theorem
\ref{kesten}.

\section{The regular case, and the proof of Theorem \ref{errw}}

We begin the section by giving some notation. Let $\P$ denote the
distribution of $\omega$ conditionally on $\t$, and $\p^x$ the law
defined by $\p^x (\cdot) := \int P_\omega^x (\cdot) \P(\! \d
\omega)$. We emphasize that $P_{\omega}^x,\,\P$ and $\p^x$ depend on
$\t$. We respectively associate the expectations $E_{\omega}^x$,
$\E$, $\e^x$. We denote also by $\Q$ and $\q^x$ the measures:
\begin{eqnarray*}
  \Q(\cdot) &:=& \int \P(\cdot) GW(\! \d \t)\,,\\
  \q^x(\cdot) &:=& \int \p^x(\cdot) GW(\! \d \t)\,.
\end{eqnarray*}

\noindent For sake of brevity, we will write $\p$ and $\q$ for
$\p^e$ and $\q^e$.

Define for $x,y\in\t$, and $n\ge 1$,
\begin{eqnarray*}
  Z_n       &:=&      \#\{x\in \t :\; |x|=n\}\,,\\
  x\le y    &\Leftrightarrow&  \exists\, p\ge 0,\; \exists\,
  x=x_0,\ldots,x_p=y\in\t\;\mbox{such that}\, \forall\, 0\le i< p\,,\; x_{i}=\buildrel
  \leftarrow \over x_{i+1}.
\end{eqnarray*}

\noindent If $x\le y$, we denote by $[\![x,y]\!]$ the set
$\{x_0,x_1,\ldots,x_p\}$, and say that $x<y$ if moreover $x\neq y$.

Define for $x\neq {\buildrel \leftarrow \over e}$, and $n\ge 1$,
\begin{eqnarray*}
  T_x         &:=&      \inf\left\{k\ge 0 :
                                            \; X_k=x \right\}\,, \\
  T_x^*       &:=&      \inf\left\{k\ge 1 : \; X_k=x
                                                  \right\}\,,\\
  \beta(x)    &:=&      P_{\omega}^x(T_{{\buildrel \leftarrow \over
                                                               x}}=\infty)\,.
\end{eqnarray*}

\noindent We observe that $\beta(x)$, $x\in \t\backslash\{{\buildrel
\leftarrow \over e}\}$, are identically distributed under $\Q$. We
denote by $\beta$ a generic random variable distributed as
$\beta(x)$. Since the walk is supposed transient, $\beta>0$
$\Q$-almost surely, and in particular $E_{\Q}[\beta]>0$.

We still consider a general Galton--Watson tree. We prove that the
number of sites visited at a generation has a bounded expectation
under $\q$.

\begin{lemma}
\label{rn} There exists a constant $c_1$ such that for any $n\ge 0$,
\begin{eqnarray*}
  E_{\q}\left[\sum_{|x|=n}\ind_{\{T_x<\infty\}}\right]\le c_1 \,.
  \label{rn2}
\end{eqnarray*}
\end{lemma}

\noindent {\it Proof.}  By the Markov property, for any $n\ge 0$,
\begin{eqnarray*}
  \sum_{|x|= n}P_{\omega}^e(T_x<\infty)\beta(x) = \sum_{|x|= n}P_{\omega}^e(T_x<\infty,~X_k\neq {\buildrel \leftarrow \over x}~\forall k> T_x)\le 1\,.
\end{eqnarray*}

\noindent The last inequality is due to the fact that there is at
most one regeneration time at the $n$-th generation. Since
$P_{\omega}^e(T_x<\infty)$ is independent of $\beta(x)$, we obtain:
$$
1 \ge E_\Q \left[ \sum_{|x|= n} P_{\omega}^e(T_x<\infty)\beta(x)
\right] = \sum_{|x|= n} E_\Q \left[ P_{\omega}^e(T_x<\infty) \right]
E_{\Q}[\beta] .
$$

\noindent In view of the identity
$E_{\q}\left[\sum_{|x|=n}\ind_{\{T_x<\infty\}}\right] = \sum_{|x|=
n} E_\Q \left[
    P_{\omega}^e(T_x<\infty) \right]$, the lemma follows immediately. $\Box$
\\ \\
Let us now deal with the case of the regular tree. We suppose in the
rest of the section that there exists $b\ge 2$ such that $\nu(x)=b$
for any $x\in\t\setminus\{{\buildrel \leftarrow\over e}\}$.
\begin{lemma}
\label{beta} If $\E\left[{1\over \sum_{i=1}^bA_i}\right]<\infty$,
then
$$
 \E\left[\frac{1}{\beta}\right]<\infty\,.
$$
\end{lemma}

\noindent {\it Proof.} Notice that $\E\left[{1\over \max_{1\le i\le
b} A_i}\right]<\infty$. For any $n\ge 0$, call $v_n$ the vertex
defined by iteration in the following way:
\begin{itemize}
\item $v_0=e$
\item $v_{n}\le v_{n+1}$ and $A(v_{n+1})=\max\{A(y),\,y\; \mbox{is a child of}\; v_n\}$.
\end{itemize}

\noindent  The Markov property tells that
$$
\beta(x)=\sum_{i=1}^{b}\omega(x,x_i)\beta(x_i)+\sum_{i=1}^{b}\omega(x,x_i)(1-\beta(x_i))\beta(x)\,
,
$$

\noindent from which it follows that for any vertex $x$,
\begin{equation}
  \label{beta2}{1 \over \beta(x)}=1+ {1\over
  \sum_{i=1}^{b}A(x_i)\beta(x_i)}\le 1+ \min_{1\le i\le b}{1\over A(x_i)\beta(x_i)}\,.
\end{equation}

\noindent Let $\mathcal{C}(v_n):=\{y\; \mbox{is a child of}\;
v_n,\,y\neq v_{n+1}\}$ be the set of children of $v_n$ different
from $v_{n+1}$. Take $C> 0$ and define for any $n\ge 1$ the event
$$E_n:=\{\forall k\in
[0,n-1]\,, \forall y\in \mathcal{C}(v_k)\,,(A(y)\beta(y))^{-1}>
C\}\,.$$

\noindent We extend the definition to $n=0$ by $E_0^c:=\emptyset$.
Notice that the sequence of events is decreasing. Using equation
(\ref{beta2}) yields
\begin{equation}
\label{beta4} {\ind_{E_{n}} \over \beta(v_n)}\le (1+ C) +
{\ind_{E_{n+1}}\over A(v_{n+1}) \beta(v_{n+1})}\,.
\end{equation}

\noindent On the other hand, by the i.i.d.\ property of the
environment, we have
$$
\P({E_n})=\P(E_1)^n\,.
$$

\noindent By choosing $C$ such that $\P(E_1)<1$ and using the
Borel--Cantelli lemma, we have $\ind_{E_n}=0$ from some $n_0\ge 0$
almost surely. Iterate equation (\ref{beta4}) to obtain
$$
\frac{1}{\beta(e)}\le \left(1 + C\right)\left(1+\sum_{n\ge
1}B(n)\right)
$$
where $B(n)=\ind_{E_n}\prod_{k= 1}^n\frac{1}{A(v_{k})}$. Hence
$$
\E\left[\frac{1}{\beta}\right]\le \left(1+C\right)\left(1+\sum_{n\ge
1}\E\left[B(n)\right]\right)\,.
$$
We observe that $
\E[B(n)]=\left\{\E\left[\ind_{E_1}A(v_1)^{-1}\right]\right\}^{n} $. When $C$ tends to infinity, $\E\left[\ind_{E_1}A(v_1)^{-1}\right]$ tends to zero since $\E[A(v_1)^{-1}]<\infty$. Choose $C$ such that $\E\left[\ind_{E_1}A(v_1)^{-1}\right]<1$ to complete the proof. $\Box$\\

For $x\in \t$ and $n\ge -1$, let
\begin{eqnarray*}
  N(x)&:=&\sum_{k\ge 0}\ind_{\{X_k=x\}}\,,\\
  N_n&:=&\sum_{|x|= n}N(x)\,,\\
  \tau_n&:=&\inf\left\{k\ge 0\, : |X_k|=n\right\}\,.
\end{eqnarray*}

\noindent In words, $N(x)$ and $N_n$ denote, respectively, the time
spent by the walk at $x$ and at the $n$-th generation, and $\tau_n$
stands for the first time the walk reaches the $n$-th generation. A
consequence of the law of large numbers is that
$$
\lim_{n\rightarrow\infty}{\tau_n\over n}= {1\over v}\;\;
\q\,\mbox{-}\,a.s.
$$

\noindent Our next result gives an upper bound for the expected
value of $N_n$.
\begin{proposition} \label{green} Suppose that $\E\left[{1\over \sum_{i=1}^b A(x_i)}\right]<\infty$. There exists a constant
$c_2$ such that for all $n\ge 0$, we have
$$
\e\left[\sum_{k=0}^nN_k\right]\le c_2\,n \,.
$$
\end{proposition}

\noindent {\it Proof.} By the strong Markov property,
$P_{\omega}^x(N(x)=\ell) = \{
P_{\omega}^x(T_x^*<\infty)\}^{\ell-1}P_{\omega}^x(T_x^*=\infty)$,
for $\ell\ge 1$. Accordingly,
$$
E_{\omega}^e\left[\sum_{k=0}^nN_k\right]=\sum_{0\le |x|\le
n}P_{\omega}^e(T_x<\infty)E_{\omega}^x[N(x)]=\sum_{0\le |x|\le
n}{P_{\omega}^e(T_x<\infty)\over {1-P_{\omega}^x(T_x^*<\infty)}} ~.
$$

\noindent We observe that $1-P_{\omega}^x(T_x^* <\infty)\,\ge\,
\sum_{i=1}^{b}\omega(x,x_i)\,\beta(x_i)$. Since
$P_{\omega}^e(T_x<\infty)$ is independent of
$(\omega(x,x_i)\beta(x_i), \; 1\le i\le b)$, we have
\begin{eqnarray}
          \e\left[\sum_{k=0}^nN_k\right]
 &\le&    \sum_{0\le |x|\le
          n} \E \left[ P_{\omega}^e(T_x<\infty)\right]\E\left[\left(\sum_{i=1}^{b}\omega(e,e_i)\beta(e_i)\right)^{\!\!-1\,}\right]
  \nonumber
  \\
 &=&      \E\left[\sum_{0\le
          |x|\le
          n}P_{\omega}^e(T_x<\infty)\right]\E\left[\left(\sum_{i=1}^{b}\omega(e,e_i)\beta(e_i)\right)^{\!\!-1\,}\right]
  \label{bet}\,.
\end{eqnarray}

\noindent Since
$\sum_{i=1}^{b}\omega(e,e_i)\,\beta(e_i)\,\ge\,\left\{\min_{i=1\ldots
b}\beta(e_i)\right\}\sum_{i=1}^{b}\omega(e,e_i) $, it follows that
\begin{eqnarray*}
  \e\left[\sum_{k=0}^nN_k\right] \le
  \E\left[\sum_{0\le |x|\le n}P_{\omega}^e(T_x<\infty)\right]
  \E\left[{1\over 1-\omega(e,{\buildrel \leftarrow \over
  e})}\right]\E\left[\left(\min_{i=1\ldots b}\beta(e_i)\right)^{\!\!-1\,}\right]\,.
\end{eqnarray*}

\noindent By definition, ${1\over 1-\omega(e,{\buildrel \leftarrow
\over e})} = 1+ {1\over \sum_{i=1}^{b} A(e_i)}$, which implies that
$\E\left[{1\over 1-\omega(e,{\buildrel \leftarrow \over
e})}\right]<\infty$. Notice also that $\E\left[\left(\min_{i=1\ldots
b}\beta(e_i)\right)^{\!\!-1\,}\right]\le b\E[{1\over \beta}]<\infty$
by Lemma \ref{beta}. Finally, use Lemma \ref{rn} to complete the
proof.
$\Box$\\

We are now able to prove the positivity of the speed.\\

\noindent {\it Proof of Theorem \ref{errw}}. We note that $\tau_n
\le \sum_{k=-1}^nN_k$ and that $N_{-1}\le N_0$. By Proposition
\ref{green}, we have $\e[\tau_n]\le 2c_2\,n$. Fatou's lemma yields
that $\e[\liminf_{n\rightarrow \infty}{\tau_n \over n}] \le 2c_2$.
Since $\lim_{n\rightarrow\infty}{\tau_n\over n}={1\over v}$, then
$v>0$. $\Box$

\section{Proof of Theorem \ref{kesten}: upper bound}

This section is devoted to the proof of the upper bound in Theorem
\ref{kesten}, which is equivalent to the following:
\begin{proposition}
\label{sup} We have
$$\liminf_{n\rightarrow\infty}\frac{\ln(\tau_n)}{\ln(n)}\ge
{1\over \Lambda} \qquad \q-a.s.$$
\end{proposition}

\subsection{Basic facts about regenerative times}

We recall some basic facts about regenerative times for the
transient RWRE. These facts can be found in \cite{gross} in the case
of regular trees, and in \cite{lpp96} in the case of biased random
walks on Galton--Watson trees.

Let
$$
D(x):=\inf\left\{k\ge 1\,: X_{k-1}=x,\, X_k={\buildrel \leftarrow
\over x} \right\},~~(\inf\emptyset:=\infty)\,.
$$

\noindent We define the first regenerative time
$$
\Gamma_1:=\inf\left\{k>0 \, : \nu(X_k)\ge 2,\,
D(X_k)=\infty,\,k=\tau_{|X_k|} \right\}
$$

\noindent as the first time when the walk reaches a generation by a
vertex having more than two children and never returns to its
parent. We define by iteration $$\Gamma_n:=\inf\left\{k>
\Gamma_{n-1} \, : \nu(X_k)\ge 2,\, D(X_k)=\infty
,\,k=\tau_{|X_k|}\right\}$$ for any $n\ge 2$ and we denote by
$\S(.)$ the conditional distribution $\q(.\,|\,\nu(e)\ge
2,\,D(e)=\infty)$.

\bigskip

\noindent {\bf Fact } {\it Assume that the walk is transient.

{\rm (i)} For any $n\ge 1$,\, $\Gamma_n<\infty$ ~~$\q$-a.s.

{\rm (ii)} Under
$\q$,\,$(\Gamma_{n+1}-\Gamma_n,|X_{\Gamma_{n+1}}|-|X_{\Gamma_n}|),\,n\ge
1$ are independent and distributed as $(\Gamma_1,|X_{\Gamma_1}|)$
under the distribution $\S$.

{\rm (iii)} We have $E_{\S}[\,|X_{\Gamma_1}|\,]<\infty$.
   }

\bigskip

We feel free to omit the proofs of {\rm (i)} and {\rm (ii)}, since
they easily follow the lines in \cite{gross} and \cite{lpp96}. To
prove {\rm (iii)}, we will show that
$E_{\S}[\,|X_{\Gamma_1}|\,]={1/E_{\Q}[\beta]}$. For any $n\ge 0$, we
have, conditionally on $|X_{\Gamma_1}|$,
\begin{eqnarray*}
    \q\left(\exists k\ge 2 \,:\, |X_{\Gamma_k}|=n\,\bigg|\, |X_{\Gamma_1}| \right)
=
    \ind_{\left\{|X_{\Gamma_1}|\le n\right\}}\q\left(\exists k\ge 2 \,:\, |X_{\Gamma_k}|-|X_{\Gamma_1}|=n-|X_{\Gamma_1}|\,\bigg|\,
    |X_{\Gamma_1}|\right).
\end{eqnarray*}

\noindent By the renewal theorem (see chapter XI of \cite{Fe} for
instance) and the fact that $\ind_{\{|X_{\Gamma_1}|\le n\}}$ tends
to $1$ $\q$-almost surely, we obtain that
\begin{eqnarray*}
\lim_{n\rightarrow\infty}\q\left(\exists k\ge 2 \,:\,
|X_{\Gamma_k}|=n\,\bigg|\, |X_{\Gamma_1}| \right) =
1/E_{\S}[\,|X_{\Gamma_1}|\,]\,.
\end{eqnarray*}

\noindent The dominated convergence yields then
\begin{eqnarray*}
\lim_{n\rightarrow\infty}\q\left(\exists k\ge 2 \,:\,
|X_{\Gamma_k}|=n \right) = 1/E_{\S}[\,|X_{\Gamma_1}|\,]\,.
\end{eqnarray*}

\noindent It remains to notice that on the other hand,
\begin{eqnarray*}
   \q\left( \exists k \in \bN \,:\,
   |X_{\Gamma_k}|=n \right)
=
   \q \left( D(X_{\tau_n})=\infty \right)
=
    E_{\Q}[\beta]\,. \;\Box
\end{eqnarray*}

If we denote for any $n\ge 0$ by $u(n)$ the unique integer such that
$\Gamma_{u(n)}\le \tau_n<\Gamma_{u(n)+1}$, then Fact yields that
$\lim_{n\rightarrow\infty}{n\over u(n)}=E_{\S}[\,|X_{\Gamma_1}|\,]$.
In turn, we deduce that
\begin{eqnarray}
\label{infeq}\liminf_{n\rightarrow\infty}\;{\ln(\tau_n)\over
\ln(n)}\; \ge \;\liminf_{n\rightarrow\infty}\;{\ln(\Gamma_n)\over
\ln(n)}\,\qquad \q\,\mbox{-}\,a.s.
\end{eqnarray}

\noindent Let for $\lambda\in [0,1]$ and $n\ge 0$,
$$
S(n,\lambda):=\sum_{k=1}^n(\Gamma_{k}-\Gamma_{k-1})^{\lambda}\,,
$$
by taking $\Gamma_0:=0$. Then $ (\Gamma_n)^{\lambda} \le
S(n,\lambda) $ since $\lambda\le 1$, which gives, by the law of
large numbers,
\begin{eqnarray}
 \label{supeq} \limsup_{n\rightarrow\infty}{(\Gamma_n)^{\lambda}\over
n}\le \lim_{n\rightarrow\infty}{S(n,\lambda)\over
n}=E_{\S}[\Gamma_1^{\lambda}]\,\qquad \q\,\mbox{-}\,a.s.
\end{eqnarray}

\subsection{Proof of Proposition \ref{sup}}

We construct a RWRE on the half-line as follows; suppose that
$\t=\{-1,0,1,\ldots\}$. This would correspond to the case where
$q_1=1$, $e=0$, ${\buildrel \leftarrow \over e}=-1$. Marking each
integer $i\ge 0$ with i.i.d.\ random variables $A(i)$, we thus
define a one-dimensional RWRE as we defined it in the case of a
Galton--Watson tree. We call $(R_n)_{n\ge 0}$ this RWRE. We still
use the notation $P^i_{\omega}$ and $\p^i$ to name the quenched and
the annealed distribution of $(R_n)$ with $R_0=i$. For $i\ge -1$ and
$a\in \mathbb{R}_+$, define $T_i:=\inf\{n\ge 0\,:\,R_n=i\}$ and
\begin{eqnarray}
\label{pna} p\,(i,a):=\p^0(T_{-1}\land T_i>a)\,,
\end{eqnarray}

\noindent where $b\land c:= \min\{ b, \, c\}$. We give two
preliminary results.

\begin{lemma}
\label{infeq2} Let $\Lambda$ be as in (\ref{Lambda}). Then
$$
\liminf_{a\rightarrow\infty}\left\{\sup_{i\ge
0}{\ln\left(q_1^ip\,(i,a)\right)\over \ln(a)}\right\}\ge -\Lambda\,.
$$
\end{lemma}

\noindent {\it Proof.} See Section 8. $\Box$ \\

We return to our general RWRE $(X_n)_{n\ge 0}$ on a general
Galton--Watson tree $\t$.
\begin{lemma}
\label{T} We have
$$
\liminf_{a\rightarrow\infty}{
\ln(\,\S\left(\Gamma_1>a\right)\,)\over \ln(a)}\ge -\Lambda\,.
$$
\end{lemma}

\noindent {\it Proof.} For any $x\in\t$, let $h(x)$ be the unique
vertex such that
\begin{eqnarray*}
x\le h(x),\qquad \nu(h(x))\ge 2\,,\qquad \forall\, y\in \t,\, x\le
y<h(x) \Rightarrow \nu(y)=1\,.
\end{eqnarray*}
In words, $h(x)$ is the oldest
 descendent of $x$ such that $\nu(h(x))\ge 2$ (and can be $x$ itself if $\nu(x)\ge 2$). We observe that
$\Gamma_1\ge T_e^*\land T_{h(X_1)}$. Moreover, $\{\nu(e)\ge
2,\,D(e)=\infty\}\supset E_1\cup E_2$ where
\begin{eqnarray*}
E_1&:=&\{\nu(e)\ge 2\}\cap\left\{X_1\neq \buildrel \leftarrow \over
e,\,T_e^*<T_{h(X_1)},\,X_{T_e^*+1}\notin\{\buildrel \leftarrow \over
e,\,X_1\}\right\}\cap\left\{X_n\neq e,\,\forall\,n\ge T_e^*+1\right\}\,,\\
E_2&:=&\{\nu(e)\ge 2\}\cap\left\{X_1\neq \buildrel \leftarrow \over
e,\,T_{h(X_1)}<T_e^*\right\}\cap\left\{X_n\neq \buildrel
\longleftarrow \over {h(X_1)},\,\forall\, n\ge
T_{h(X_1)}+1\right\}\,.
\end{eqnarray*}

\noindent It follows that
\begin{eqnarray}
  \S(\Gamma_1>a)\ge {1\over \q(\nu(e)\ge 2,\,D(e)=\infty)}(
\q(T_e^*>a,\,E_1)+\q(T_{h(X_1)}>a,\,E_2))\label{T^*}\,.
\end{eqnarray}

\noindent We claim that
\begin{eqnarray}
\label{eq1} \q(T_e^*>a,\,E_1)= c_{3}\,\q(T_{\buildrel
\leftarrow\over e}< T_{h(e)},\,1+T_{\buildrel \leftarrow\over
e}>a)\,.
\end{eqnarray}

\noindent Indeed, write
\begin{eqnarray*}
 P_{\omega}^e(T_e^*>a,\,E_1)  = \sum_{e_i\neq
  e_j}P_{\omega}^e\left(T_e^*<T_{h(e_i)},\,X_1=e_i,\,X_{T_e^*+1}=e_j,\,D(e_j)=\infty,\,T_e^*>a\right).
\end{eqnarray*}

\noindent By gradually applying the strong Markov property at times
$T_e^*+1$, $T_e^*$ and at time $1$, this yields
\begin{eqnarray*}
   P_{\omega}^e(T_e^*>a,\,E_1)=
   \sum_{e_i\neq
   e_j}\omega(e,e_i)P_{\omega}^{e_i}\left(T_e<T_{h(e_i)},\,1+T_e>a\right)\omega(e,e_j)\beta(e_j).
\end{eqnarray*}

\noindent Since $\omega(e,e_i)\omega(e_j)$, $\beta(e_j)$ and
$P_{\omega}^{e_i}\left(T_e<T_{h(e_i)},\,1+T_e>a\right)$ are
independent under $\P$, this leads to
\begin{eqnarray*}
  \p(T_e^*>a,\,E_1)=
  \sum_{e_i\neq
   e_j}\E\left[\omega(e,e_i)\omega(e,e_j)\right]\p^{e_i}\left(T_{e}<T_{h(e_i)},\,1+T_e>a\right)\E\left[\beta(e_j)\right].
\end{eqnarray*}

\noindent By the Galton--Watson property,
\begin{eqnarray*}
  \q(T_e^*>a,\,E_1)=
  E_{\Q}\left[\ind_{\{\nu(e)\ge 2\}}\sum_{e_i\neq
   e_j}\omega(e,e_i)\omega(e,e_j)\right]\q^{e}\left(T_{\buildrel \leftarrow \over
   e}<T_{h(e)},\,1+T_{\buildrel \leftarrow \over e}>a\right)E_{\Q}\left[\beta\right]\,,
\end{eqnarray*}

\noindent which gives (\ref{eq1}). Similarly,
\begin{eqnarray}
\q(T_{h(X_1)}>a,\,E_2)=c_{4}\q\left(T_{\buildrel \leftarrow\over
e}>T_{h(e)},\,1+T_{h(e)}>a\right)\label{t2}\,.
\end{eqnarray}

\noindent Finally, by (\ref{T^*}), (\ref{eq1}) and (\ref{t2}) we get
\begin{eqnarray*}
\S(\Gamma_1>a) \ge c_{5}\,\q\left(1+T_{\buildrel \leftarrow \over
e}\land T_{h(e)}>a\right).
\end{eqnarray*}

\noindent Conditionally on $|\,h(e)|$, the walk $|\,X_n|,\,0\le n\le
T_{\buildrel \leftarrow \over e}\land T_{h(e)} $ has the
distribution of the walk $R_n,\, 0\le n\le T_{-1}\land T_{|h(e)|}$,
as defined at the beginning of this section. For any $n\ge 0$, since
$GW(|h(e)|=n)=q_1^{n}(1-q_1)$, it follows that
$\q\left(1+T_{\buildrel \leftarrow \over e}\land
T_{h(e)}>a\right)\ge q_1^n(1-q_1)p\,(n,a)$. Finally,
\begin{eqnarray*}
\liminf_{a\rightarrow\infty}{\ln\left(\S\left(\Gamma_1>a\right)\right)\over
\ln(a)}\ge \liminf_{a\rightarrow\infty}\left\{\sup_{n\ge
0}{\ln\left(q_1^np\,(n,a)\right)\over \ln(a)}\right\}\,.
\end{eqnarray*}
Applying Lemma \ref{infeq2} completes the proof.  $\Box$

\bigskip

We now have all of the ingredients needed for the proof of
Proposition \ref{sup}.

\bigskip

\noindent {\it Proof of Proposition \ref{sup}.} If $\Lambda \ge 1$,
Proposition \ref{sup} trivially holds since $\tau_n\ge n$. We
suppose that $\Lambda<1$, and let $\Lambda<\lambda<1$. Let
$M_n:=\max\{\Gamma_k-\Gamma_{k-1},\;k=2,\ldots n\}$. We have
$\q\left(M_n\le n^{1\over
\lambda}\right)=\q\left(\Gamma_2-\Gamma_1\le n^{1\over
\lambda}\right)^n$. By Lemma \ref{T}, $\q\left(\Gamma_2-\Gamma_1\le
n^{1\over \lambda}\right)\le 1-n^{-1+\varepsilon}$ for some
$\varepsilon>0$ and large $n$. Consequently, $\sum_{n\ge
1}\q\left(M_n\le n^{1\over \lambda}\right)<\infty$, and the
Borel-Cantelli lemma tells that $\q$-almost surely and for
sufficiently large $n$, $M_n\ge n^{1\over \lambda}$, which in turn
implies that $\liminf_{n\rightarrow\infty}{\Gamma_n-\Gamma_1\over
n^{1\over \lambda}}\ge 1$. We proved then that
$\liminf_{n\rightarrow\infty}{\ln(\Gamma_n)\over \ln(n)}\ge {1\over
\Lambda}$. Therefore, by equation (\ref{infeq}),
$$\liminf_{n\rightarrow\infty}{\ln(\tau_n)\over \ln(n)}\ge
{1\over \Lambda}\,\qquad \q\,\mbox{-}\,a.s.\; \Box$$

\section{Technical results}

We give, in this section, some tools needed in our proof of the
lower bound in Theorem \ref{kesten}. $Z_n$ stands as before for the
size of the $n$-th generation of $\t$.

\begin{lemma}
\label{growth} For every $b,n\ge 1$, we have
$$
E_{GW}[Z_n\ind_{\{Z_n\le b\}}]\le bn^bq_1^{n-b}\,.
$$
\end{lemma}

\noindent {\it Proof.} If $Z_n\le b$, then there are at most $b$
vertices before the $n$-th generation having more than one child.
Therefore,
$$
GW(Z_n\le b)\le C_n^bq_1^{n-b}\le n^bq_1^{n-b}
$$
and we conclude since $E_{GW}[Z_n\ind_{\{Z_n\le b\}}]\le b\,GW(Z_n\le b)$. $\Box$\\

\begin{lemma} \label{b}  Let $\beta_i$, $i\ge 1$ be independent random variables distributed as $\beta$. There exists $b_0\ge 1$ such that
$$
E_{\Q}\left[\left(\frac{1}{\sum_{i=1}^{b_0}\beta_i}\right)^{\!\!2}\right]<\infty\,.
$$
\end{lemma}

\noindent {\it Proof.} Let $\t^{(i)}$, $i\ge 1$ be independent
Galton--Watson trees of distribution $GW$. We equip independently
each $\t^{(i)}$ with an environment of distribution $\P$ so that we
can look at the random variable $\beta(e^{(i)})$ where $e^{(i)}$ is
the root of $\t^{(i)}$. Then $\beta(e^{(i)})$, $i\ge 1$ are
independent random variables distributed as $\beta$.

Let $c_{6}> 0$ be such that $\eta := \Q(\frac{1}{\beta}> c_{6})<1$.
Recall that $\frac{1}{\alpha}\le A(x)\le \alpha$, $\forall x\in \t$,
$\q$-almost surely. Let $R^{(i)}:=\inf\{n\ge 0 \,: \exists
y\in\t^{(i)},\,|y|=n, \frac{1}{\beta(y)}\le c_{6}\}$ be the first
generation in $\t^{(i)}$ where a vertex verifies
$\frac{1}{\beta(y)}\le c_{6}$, and let $y^{(i)}$ be such a vertex
$y$. Recall from equation (\ref{beta2}) that
$$
{1\over \beta(x)}\le 1+ {1\over A(x_j)\beta(x_j)}
$$
for any child $x_j$ of a vertex $x$. By iterating the inequality
 on the path $[\![e^{(i)},y^{(i)}]\!]$, we obtain
$$
\frac{1}{\beta(e^{(i)})}\le 1 + \sum_{z\in ]\!]e,y^{(i)}[\![}H(z) +
\frac{H(y^{(i)})}{\beta(y^{(i)})}
$$
where $H(z)=\prod_{v\in ]\!]e^{(i)},z]\!]}\frac{1}{A(v)}\le
\alpha^{|z|}$  for every $z\in\t$ by the bound assumption on $A$.
Since ${1\over \beta(y^{(i)})}\le c_{6}$, this implies
\begin{eqnarray*}
  \frac{1}{\beta(e^{(i)})}\le c_{7}\,\alpha^{R^{(i)}}\, ,
  \label{beta1}
\end{eqnarray*}

\noindent for some constant $c_{7}$. There exist constants $c_8$ and
$c_9$ such that for any $b\ge 1$,
\begin{eqnarray}
     \label{beta5} \left(\frac{1}{\sum_{i=1}^{b}\beta(e^{(i)})}\right)^{\!\!2}\le c_8\,c_9^{\min_{1\le i\le b}R^{(i)}}\,.
\end{eqnarray}

\noindent We observe that
\begin{eqnarray}
     E_{\Q}\left[c_9^{\min_{1\le i\le b}R^{(i)}}\right]&=&\sum_{n=0}^{\infty}c_9^n\,\Q(\min_{1\le i\le b}R^{(i)}=n)\nonumber \\
&\le& \sum_{n=0}^{\infty}c_9^n\,\Q(R^{(1)}\ge n)^b\,.\label{sum}
\end{eqnarray}

\noindent We have, for any $n\ge 1$, $ \Q(R^{(1)}\ge n)   \le
\Q\left(\forall |x|=n-1, \frac{1}{\beta(x)}>c_{6}\right) $. Recall
that $\eta :=\Q(\frac{1}{\beta}> c_{6})<1$. By independence,
$$
\Q\left(\forall |x|=n-1,
\frac{1}{\beta(x)}>c_6\right)=E_{GW}[\eta^{Z_{n-1}}]\,.
$$

\noindent Let $q_1<a<1$. There exists a constant $c_{10}$ such that
$E_{GW}\left[\eta^{Z_{\ell}}\right]\le c_{10}\,a^{\ell+1}$ for any
$\ell\ge 0$. Choose $b_0$ such that $c_9a^{b_0}<1$. Then by
(\ref{sum}), $E_{\q}\left[c_9^{\min_{1\le i\le
b_0}R^{(i)}}\right]<\infty$,
which completes the proof in view of (\ref{beta5}). $\Box$ \\

Define for any $u,v\in\t$ such that $u\le v$ and for any $n\ge 1$:
\begin{eqnarray}
    \label{defp1}
    p_1(u,v)
&=& P_{\omega}^u\left(
    T_{{\buildrel \leftarrow \over u}}=\infty\,,\,T_u^*=\infty\,,\,T_v=\infty
    \right)\,,\\ \label{nu}
    \nu(u,n)
&=& \#\left\{x\in\t\,:\,u\le x, |\,x-u|=n \right\}\,.
\end{eqnarray}

\begin{lemma}
\label{p1}  For all $n\ge 2$ and $k\in\{1,2\}$, we have
\begin{eqnarray}
E_{\Q}\left[\sum_{|u|=n}\frac{\ind_{\{Z_n>b_0\}}}{[p_1(e,u)]^k}\right]
< \infty\, \label{p11}.
\end{eqnarray}
\end{lemma}

\noindent {\it Proof.} Let $n\ge 2$ and $k\in\{1,2\}$ be fixed
integers and $\widetilde{n}:=\inf\{\ell\ge 1\;:\;Z_{\ell}>b_0\}$.
Notice that $\{Z_n>b_0\}=\{\widetilde n\le n\}$. For any $u\in\t$
such that $|u|\ge \widetilde n$, let $\widetilde u \in \t$ be the
unique vertex such that $|\widetilde u|=\widetilde n$ and
$\widetilde u\le u$ that is the ancestor of $u$ at generation
$\widetilde n$. We have by the Markov property,
\begin{eqnarray}
   p_1(e,u)
\ge
   \sum_{|y|=\widetilde{n}-1}P_{\omega}^e(T_y<T_e*)P_{\omega}^y(T_{{\buildrel \leftarrow \over
   y}}=\infty\,,\,T_{\widetilde u}=\infty). \label{p1e}
\end{eqnarray}

\noindent For any $|y|\le \widetilde{n}$ and $y_i$ child of $y$, we
observe that
\begin{eqnarray*}
\label{omega} \omega(y,y_i)={A(y_i) \over 1+\sum_{j=1}^{\nu(y)}
A(y_j)}\ge {1\over c_{11}\nu(y)}\,, \end{eqnarray*}

\noindent which is greater than  $1/c_{11}b_0:=c_{12}$, by the
boundedness assumption on $A$ and the definition of $\widetilde{n}$.
It yields that for any $|y|=\widetilde{n}-1$,
\begin{equation}
P_{\omega}^e(T_y<T_e^*)\ge P_{\omega}^e(X_{\widetilde n -1}=y)\ge
c_{12}^{\widetilde n}\,. \label{p1a}
\end{equation}

\noindent By the Markov property,
\begin{eqnarray*}
  &&P_{\omega}^y(T_{\buildrel \leftarrow \over
  y}=\infty,T_{y_i}=\infty)\\
&=&
  \sum_{j\neq
  i}\omega(y,y_j)\beta(y_j)+\left(\sum_{j\neq
  i}\omega(y,y_j)(1-\beta(y_j))\right)P_{\omega}^y(T_{\buildrel \leftarrow
  \over y}=\infty,T_{y_i}=\infty).
\end{eqnarray*}

\noindent This leads to
\begin{eqnarray*}
   P_{\omega}^y(T_{\buildrel \leftarrow \over
   y}=\infty,T_{y_i}=\infty)
&=&
   \frac{\sum_{j\neq i} A(y_j)\beta(y_j)}{1+A(y_i)+\sum_{j\neq i}
   A(y_j)\beta(y_j)} \nonumber\\
&\ge&
   {1\over \alpha(1+\alpha)}{\frac{\sum_{j\neq i}\beta(y_j)}{1+\sum_{j\neq i}
   \beta(y_j)}}\\
&\ge&
   {1 \over 2\alpha(1+\alpha)}\left(1\land \sum_{j\neq i}\beta(y_j)\right)\,.
\end{eqnarray*}

\noindent Similarly, $P_{\omega}^y(T_{\buildrel \leftarrow \over
   y}=\infty)\ge {1\over 2\alpha^2}{\left(1\land
   \sum_{j=1}^{\nu(y)}\beta(y_j)\right)}$. Thus, we have for any $|y|=\widetilde
   n -1$,
\begin{eqnarray}
\label{p1b} P_{\omega}^y(T_{{\buildrel \leftarrow \over
   y}}=\infty\,,\,T_{\widetilde u}=\infty)\ge c_{13}\left(1\land \sum_{y_j\neq \widetilde u}\beta(y_j)\right)\,.
\end{eqnarray}

\noindent By equations (\ref{p1e}), (\ref{p1a}) and (\ref{p1b}), we
have
\begin{eqnarray*} \label{p1d} p_1(e,u)\ge
c_{13}c_{12}^{\widetilde n}\left(1\land \sum_{|x|=\widetilde n
:x\neq \widetilde u}\beta(x)\right).
\end{eqnarray*}

\noindent Therefore, arguing over the value of $\widetilde u$, we
obtain
\begin{eqnarray*}
\label{p1c} \ind_{\{n\ge \widetilde
n\}}\sum_{|u|=n}\E\left[\frac{1}{[p_1(e,u)]^k}\right] \le
c_{14}\sum_{|y|=\widetilde n}\nu(y,n-\tilde n)\E\left[1\lor {1\over
[\sum_{|x|=\widetilde{n},x\neq y}\beta(x)]^k}\right],
\end{eqnarray*}

\noindent where $c_{14}:=(c_{13}c_{12}^{n})^{-k}$. By using the
Galton--Watson property at generation $\widetilde n$,
\begin{eqnarray*}
&&
     \sum_{|u|=n}E_{\Q}\left[\frac{\ind_{\{u\in\t,Z_n>b_0\}}}{[p_1(e,u)]^k}\,\bigg|\,\widetilde
     n\,,Z_0,\ldots,Z_{\widetilde n}\right]\\
&\le&
     c_{14}\sum_{|y|=\widetilde n}E_{GW}[\nu(y,n-\widetilde
     n)]E_{\Q}\left[1 \lor {1\over
     [\sum_{i=1}^{p}\beta(i)]^k}\right]_{p=Z_{\widetilde
     n}-1}\\
&\le&
     c_{15}Z_{\widetilde n}
\end{eqnarray*}
by Lemma \ref{b}. Integrating over $GW$ completes the proof of
(\ref{p11}). $\Box$ \\

\noindent {\bf Remark.} Lemma \ref{p1} tells in particular that
\begin{eqnarray}
\label{fin} E_{\Q}\left[{\ind_{\{Z_{n}>b_0\}}\over
\beta(e)}\right]\le E_{\Q}\left[{\ind_{\{Z_{n}>b_0\}}\over
P_{\omega}^e(T_{\buildrel \leftarrow \over e}=\infty,T_e^*=\infty)}
\right]<\infty\,.
\end{eqnarray}

We deal now with a comparison between RWREs on a tree and
one-dimensional RWREs already used in \cite{lpp96}. Let $\t$ be a
tree and $\omega$ the environment on this tree. Take $x\le y\in\t$.
We look at the path $[\![{\buildrel \leftarrow \over
x},y]\!]=\{{\buildrel \leftarrow \over x}=x_{-1},x_0,\ldots,x_p=y\}$
defined as the shortest path from ${\buildrel \leftarrow \over x}$
to $y$, and we consider on it the random walk $(\widetilde{X}_n)$
with probability transitions $\widetilde{\omega}({\buildrel
\leftarrow \over x},x)=\widetilde{\omega}(y,x_{p-1})=1$ and for any
$0\le i<p$,
\begin{eqnarray*}
\widetilde{\omega}(x_i,x_{i+1})&=&{\omega(x_i,x_{i+1})\over \omega(x_i,x_{i+1})+\omega(x_i,x_i-1)}\,,\\
\widetilde{\omega}(x_i,x_{i-1})&=&{\omega(x_i,x_{i-1})\over
\omega(x_i,x_{i+1})+\omega(x_i,x_{i-1})}\,.
\end{eqnarray*}

\noindent Thus we can associate to the pair $(x,y)$ a
one-dimensional RWRE on $[\![{\buildrel \leftarrow \over x},y]\!]$,
and we denote by $\widetilde{P},\,\widetilde{E}$ the probabilities
and expectations related to this new RWRE. We observe that under
$\q^x$, the RWRE $(\widetilde{X}_n,\,n\le T_{{\buildrel \leftarrow
\over x}}\land T_{y})$ has the distribution of the RWRE $(R_n,\,n\le
T_{-1}\land T_{p})$ introduced in Section 3.2. For any $x,y\in \t$,
the event $\{T_x<T_y\}$
 means that $T_x<\infty$ and $T_x<T_y$.

\begin{lemma}\label{one} For any $x,y,z\in\t$ with $x\le z< y$,
\begin{eqnarray*}
    P_{\omega}^z(T_y<T_{\buildrel \leftarrow \over x})   &\le&   \widetilde{P}_{\omega}^z(T_y<T_{\buildrel \leftarrow \over x})\,,\\
    P_{\omega}^z(T_{\buildrel \leftarrow \over x}<T_y)   &\le&   \widetilde{P}_{\omega}^z(T_{\buildrel \leftarrow \over x}<T_y)\,.
\end{eqnarray*}
\end{lemma}

\noindent {\it Proof.} Fix $z_1,\ldots z_{n-1}$ in $]\!]{\buildrel
\leftarrow \over x},y[\![$ and $z_n\in[\![{\buildrel \leftarrow
\over x},y]\!]$. Then
\begin{eqnarray*}
     P_{\omega}^z(X_1=z_1,\ldots,X_n=z_n)
=\frac{\omega(z,z_1)}{1-f(z)}\ldots\frac{\omega(z_{n-1},z_n)}{1-f(z_{n-1})}
\end{eqnarray*}
where $f(r)$ represents the probability of making an excursion away
from the path $[\![{\buildrel \leftarrow \over x},y]\!]$ from the
vertex $r$. For each $r\in [\![{\buildrel \leftarrow \over
x},y[\![$, call $r^+$ the child of $r$ which lies in the path. Then
$f(r)\le 1-\omega(r,r^+)-\omega(r,{\buildrel \leftarrow\over r})$.
It follows that
\begin{eqnarray*}
      P_{\omega}^z(X_1=z_1,\ldots,X_n=z_n)
&\le& \widetilde{\omega}(z,z_1)
      \ldots   \widetilde{\omega}(z_{n-1},z_n) \\
&=& \widetilde{P}_{\omega}^z (
\widetilde{X}_1=z_1,\ldots,\widetilde{X}_n=z_n )\,.
\end{eqnarray*}

\noindent It remains to see that  the events $\{T_y<T_x\}$ and
$\{T_x<T_y\}$ can be written as an union of disjoint sets of the
form $\{X_1=z_1,\ldots,X_n=z_n\}$. $\Box$\\

The last lemma deals with the one-dimensional RWRE $(R_n)_{n\ge 0}$
defined in Section 3.2.

\begin{lemma}
\label{estimation} For any $n\ge 1$, there exists a number
$c_{19}(n)$
 such that for any $i> n$ and almost every
$\omega$,
$$
E_{\omega}^{0}[T_{-1}\land T_{i}]\le c_{19}
E_{\omega}^{n}[T_{n-1}\land T_{i}]\,.
$$
\end{lemma}

\noindent {\it Proof.} Let $i> n\ge 1$. By the Markov property and
for $0<p\le i$, we have
\begin{eqnarray*}
    E_{\omega}^{p-1}[T_{p-2}\land T_{i}]
=   1 + \omega (p-1,p) \left\{ E_{\omega}^{p} [ T_{p-1} \land T_{i}
    ] + P_{\omega}^{p} ( T_{p-1} < T_{i} )E_{\omega}^{p-1}[ T_{p-2} \land T_{i} ] \right\}
\end{eqnarray*}
which gives that $E_{\omega}^{p-1}[T_{p-2}\land
T_i]={1+\omega(p-1,p)E_{\omega}^p[T_{p-1}\land T_i]\over
1-\omega(p-1,p)P_{\omega}^p(T_{p-1}\land T_i)}$, so that for some
$c_{20}, c_{21}$ and $c_{22}$ we have
$$
E_{\omega}^{p-1}[T_{p-2}\land T_{i}]\le c_{20} + c_{21}
E_{\omega}^{p}[T_{p-1}\land T_{i}]\le
c_{22}E_{\omega}^{p}[T_{p-1}\land T_{i}].
$$
Iterating the inequality over all $p$ from $1$ to $n$ gives the desired inequality. $\Box$\\

\section{Proof of Theorem \ref{kesten}: lower bound}

Let $(R_n)_{n\ge 0}$ be the one-dimensional RWRE associated with
$\t=\{-1,0,1,\ldots\}$ defined in Section 3.2 and $T_{i}=\inf\{k\ge
0:\, R_k=i\}$. Define for any $\lambda\in [0,1]$,
\begin{eqnarray}
\label{mn} m(n,\lambda):=\E\left[\left(E_{\omega}^0\left[T_{-1}\land
T_n\right]\right)^{\lambda}\right]\,,
\end{eqnarray}

\noindent and let
\begin{equation}
  \label{qc} \lambda_c:=\sup\left\{\lambda\ge 0\,:\,\exists r>q_1\; \mbox{such that}\; \sum_{n\ge
  0}m(n,\lambda)r^n<\infty\right\}\,.
\end{equation}

\noindent We start with a lemma.

\begin{lemma}
\label{tec1} We have $\Lambda\le \lambda_c\,$.
\end{lemma}

\noindent {\it Proof.} See Section 8. $\Box$

\bigskip

Take a $\lambda\in [0,1]$ such that $\lambda<\Lambda$. By Lemma
\ref{tec1}, we have $\lambda<\lambda_c$ which in turn implies by
(\ref{qc}) that there exists an $1>r>q_1$ such that
\begin{equation}
     \label{lambda} \sum_{n\ge 0}m(n,\lambda)\,(n+1)r^n<\infty\,.
\end{equation}

\noindent Recall the definition of $b_0$ in Lemma \ref{b}. Then, by
Lemma \ref{growth}, we can define
\begin{eqnarray*}
     \label{defn0}
     n_0
:=   \inf\left\{n\ge 1\,:\, E_{GW}[Z_{n}\ind_{\{Z_{n}\le b_0\}}]\le
     r^{n} \right\}\,.
\end{eqnarray*}

\noindent Let $\t_{n_0}$ be the subtree of $\t$ defined as follows:
$y$ is a child of $x$ in $\t_{n_0}$ if $x\le y$ and $|y-x|=n_0$. In
this new Galton--Watson tree $\t_{n_0}$, we define
\begin{equation}
\mathbb{W}=\mathbb{W}(\t):=\{x\in\t_{n_0}: \forall y\in\t_{n_0},
(y<x) \Rightarrow \nu(y,n_0)\le b_0 \}\,, \label{defw}
\end{equation}

\noindent where $\nu(y,n_0)$ is defined in (\ref{nu}). We call $W_k$
the size of the $k$-th generation of $\mathbb{W}$. The subtree
$\mathbb{W}$ is a Galton--Watson tree, whose offspring distribution
is of mean $E_{GW}[Z_{n_0}\ind_{\{Z_{n_0}\le b_0\}}]\le r^{n_0}$. In
particular, we have for any $k\ge 0$,
\begin{eqnarray}
    \label{W} E_{GW}[W_k]\le r^{kn_0}\,.
\end{eqnarray}

\noindent For any $y\in\t$, we denote by $y_{n_0}$ the youngest
ancestor of $y$ belonging to $\t_{n_0}$, or equivalently the unique
vertex such that
\begin{eqnarray*}
y_{n_0}\le y,\qquad y_{n_0}\in \t_{n_0},\qquad \forall\;
z\in\t_{n_0}\;\; z\le y\Rightarrow z\le y_{n_0}\,.
\end{eqnarray*}

\noindent Let
\begin{eqnarray*}
   N_{1,n}&:=&\sum_{|y|=n}N(y)\ind_{\{\nu(y_{n_0},n_0)> b_0\}}\,,\\
   N_{2,n}&:=&\sum_{|y|=n}N(y)\ind_{\{\nu(y_{n_0},n_0)\le b_0,y_{n_0}\notin
   \mathbb{W}\}}\,.
\end{eqnarray*}

\begin{lemma}
\label{lem} There exists a constant $L$ such that for any $n\ge n_0$
:
\begin{eqnarray}
   E_{\q}[N_{1,n}]&\le& L \label{N1}\,,\\
   E_{\q}[N_{2,n}^{\lambda}]&\le& L \label{N2}\,.
\end{eqnarray}
\end{lemma}

We admit Lemma \ref{lem} for the time being, and show how it implies Theorem \ref{kesten}.\\

\noindent {\it Proof of  Theorem \ref{kesten}: lower bound}. Notice
that $\mathbb{W}$ is finite almost surely. Then, there exists a
random $K\ge 0$ such that for $n\ge K$, $N_n\le N_{1,n}+N_{2,n}$.
Lemma \ref{lem} yields that $E_{\q}[N_{n}^{\lambda}\ind_{\{n\ge
K\}}]\le L^{\lambda}+L$ for any $n\ge n_0$. By Fatou's lemma,
$\liminf_{n\rightarrow \infty}{\sum_{k=K}^{n} N_k^{\lambda} \over
n}<\infty$. Denote by $(r_k,\,k\ge 0)$ the sequence
$(|X_{\Gamma_{k}}|,\,k\ge 0)$. Notice that for any $k\ge 1$,
$$\Gamma_{k+1}-\Gamma_k = \sum_{i=r_k+1}^{r_{k+1}}N_i\,.$$

\noindent It yields that
$S(u(n),\lambda):=\sum_{k=1}^{u(n)}(\Gamma_{k}-\Gamma_{k-1})^{\lambda}\le
\sum_{i=0}^{r_{u(n)}}N_i^{\lambda}\le \sum_{i=0}^{n}N_i^{\lambda}$
where, as in Section 3, $u(n)$ is the unique integer such that
$\Gamma_{u(n)}\le \tau_n<\Gamma_{u(n)+1}$. Observe also that
${n\over u(n)}$ tends to
$E_{\mathbb{S}}\left[|X_{\Gamma_1}|\right]$. It follows that
$$\liminf_{n\rightarrow\infty} {1\over n}S(n,\lambda)\le
\liminf_{n\rightarrow\infty}{1\over u(n)} \sum_{k=K}^{n}
N_k^{\lambda} =
E_{\mathbb{S}}\left[|X_{\Gamma_1}|\right]\liminf_{n\rightarrow\infty}{1\over
n} \sum_{k=K}^{n} N_k^{\lambda} <\infty\,.$$

\noindent Using equation (\ref{supeq}) implies that
$\limsup_{n\rightarrow\infty}{(\Gamma_n)^{\lambda}\over n}<c_{23}$
for some constant $c_{23}$. We check that $|X_n|\ge
\#\{k\,:\,\Gamma_k\le n\}$ which leads to $|X_n|\ge
{n^{\lambda}\over c_{23}}$ for sufficiently large $n$. Letting
$\lambda$ go to $\Lambda$ completes the proof. $\Box$

\bigskip

The rest of this section is devoted to the proof of Lemma \ref{lem}.
For the sake of clarity, the two estimates, (\ref{N1}) and
(\ref{N2}), are proved in distinct parts.

\subsection{Proof of Lemma \ref{lem}: equation (\ref{N1})}

For all $y\in\t$, call $Y$ the youngest ancestor of $y$ such that
$\nu(Y,n_0)>b_0$. We have
$$
E_{\omega}^e[N(y)]=P_{\omega}^e(T_y<\infty)E_{\omega}^y[N(y)]\le
P_{\omega}^e(T_{Y}<\infty)E_{\omega}^y[N(y)]\,.
$$

\noindent We compute $E_{\omega}^y[N(y)]$ with a method similar to
the one given in \cite{lpp96}. By the Markov property,
$$
E_{\omega}^y[N(y)]= G(y,Y) + P_{\omega}^y(T_{
Y}<\infty)P_{\omega}^{Y}(T_y<\infty)E_{\omega}^y[N(y)] \,,
$$

\noindent where
$G(y,Y):=E_{\omega}^y\left[\sum_{k=0}^{T_{Y}}\ind_{\{X_k=y\}}\right]$.
When $\nu(y_{n_0},n_0)>b_0$, there exists a constant $c_{24}>0$ such
that $P_{\omega}^y(T_y^*>T_{Y})\ge c_{24}$. Therefore, in this case
$G(y,Y)\le (c_{24})^{-1}=:c_{25}$. It follows that
\begin{eqnarray*}
E_{\omega}^y[N(y)]\ind_{\{\nu(y_{n_0},n_0)>b_0\}} &\le&
c_{25}{\ind_{\{\nu(y_{n_0},n_0)>b_0\}}\over
1-P_{\omega}^{Y}(T_y<\infty)P_{\omega}^y(T_{Y}<\infty)}\\
&\le&  c_{25}{\ind_{\{\nu(y_{n_0},n_0)>b_0\}}\over
1-P_{\omega}^{Y}(T_{Y}^*<\infty)}\\
&\le& c_{25}{\ind_{\{\nu(y_{n_0},n_0)>b_0\}}\over \gamma(Y)}\,,
\end{eqnarray*}

\noindent where $\gamma(x):=P_{\omega}^x(T_{\buildrel \leftarrow
\over x}=\infty,T_x^*=\infty)$. Arguing over the value of $Y$ yields
that
\begin{eqnarray*}
E_{\q}[N_{1,n}]&\le&  c_{25}E_{\Q}\left[\sum_{n-n_0<
|z|\le n}P_{\omega}^e(T_{z}<\infty){\ind_{\{\nu(z,n_0)>b_0\}}\over \gamma(z)}\right]\\
&=&  c_{25}E_{\Q}\left[\sum_{n-n_0< |z|\le
n}P_{\omega}^e(T_{z}<\infty)\right]E_{\Q}\left[{\ind_{\{Z_{n_0}>b_0\}}\over
\gamma(e)}\right]\\
&\le& c_{25}n_0\,c_1\,c_{26}\,,
\end{eqnarray*}

\noindent by Lemma \ref{rn} and equation (\ref{fin}). $\Box$

\subsection{Proof of Lemma \ref{lem}: equation (\ref{N2})}

\noindent For any $y\in \t$ such that $\nu(y_{n_0},n_0)\le b_0$ and
$y_{n_0}\notin \mathbb{W}$, choose $Y_1=Y_1(y)$, $Y_2=Y_2(y)$ and
$Y_3=Y_3(y)$,
 vertices of $\t_{n_0}$,  such that
\begin{eqnarray*}
Y_1< y,\;\;\;\; &\nu(Y_1,n_0)> b_0,&\;\;\;\; \forall\;z\in\t_{n_0},\;Y_1< z\le y \Rightarrow \nu(z,n_0)\le b_0 \\
Y_1< Y_2\le y,\;\;\;\; &&\;\;\;\;\forall\; z\in\t_{n_0},\; Y_1<z\le y\Rightarrow Y_2\le z\,,\\
y\le Y_3,\;\;\;\; &\nu(Y_3,n_0)> b_0,&\;\;\;\; \forall\;
z\in\t_{n_0},\;y\le z<Y_3\Rightarrow \nu(z,n_0)\le b_0\,.
\end{eqnarray*}

\noindent By definition, $Y_1$ is the youngest ancestor of $y$ in
$\t_{n_0}$ such that $\nu(Y_1,n_0)>b_0$ and $Y_2$ the child of $Y_1$
in $\t_{n_0}$ which is also an ancestor of $y$. In the rest of the
section, $\widetilde P_{\omega}=\widetilde{P}_{\omega}(Y_1,Y_3)$ and
$\widetilde E_{\omega}=\widetilde{E}_{\omega}(Y_1,Y_3)$ represent
the probability and expectation for the one-dimensional RWRE
associated to the path $[\![Y_1,Y_3]\!]$, as seen in Lemma
\ref{one}. They depend then on the pair $(Y_1,Y_3)$, which doesn't
appear in the notation for sake of brevity. Define for any $n\ge
n_0$,
\begin{eqnarray}
\label{defs}
   S(n):=E_{\Q}\left[\sum_{|y|=n:Y_1=e}{\left[p_1(e,Y_2)^2\beta(Y_3)\right]^{-1}} \left(\widetilde{E}^{Y_2}_{\omega}[T_{{\buildrel \leftarrow\over Y_2}}\land
   T_{Y_3}]\right)^{\lambda}\right]\,,
\end{eqnarray}
where ${\buildrel \leftarrow\over Y_2}$ represents as usual the
parent of $Y_2$ in the tree $\t$ and $p_1(u,v)$ is defined in
(\ref{defp1}).

\begin{lemma}
\label{N31} There exists a constant $c_{27}$ such that for any $n\ge
n_0$,
\begin{eqnarray*}
    E_{\q}[N_{2,n}^{\lambda}]
\le c_{27}\sum_{k\ge n_0}S(k)\,.
\end{eqnarray*}
\end{lemma}

\noindent {\it Proof.} We observe that
\begin{eqnarray*}
E_{\omega}^e[N_n^{\lambda}]=E_{\omega}^e\left[\left(\sum_{|y|=n}N(y)\right)^{\lambda}\right]
\le E_{\omega}^e\left[\sum_{|y|=n}N(y)^{\lambda}\right]
\end{eqnarray*}
since $\lambda\le 1$. By the Markov property,
$E_{\omega}^e[\sum_{|y|=n}N(y)^{\lambda}]=\sum_{|y|=n}P_{\omega}^e(T_y<\infty)E_{\omega}^y[N(y)^{\lambda}]$.
An application of Jensen's inequality yields that
\begin{eqnarray}
\label{sumgen} E_{\omega}^e[N_n^{\lambda}] \le
\sum_{|y|=n}P_{\omega}^e(T_y<\infty)\left(E_{\omega}^y[N(y)]\right)^{\lambda}\,.
\end{eqnarray}

\noindent Using the Markov property for any $|y|=n$, we get
\begin{eqnarray*}
&&   E_{\omega}^y[N(y)]\\
&=&
     G(y,Y_1\land Y_3) + E_{\omega}^y[ N(y) ] ( P_{\omega}^y( T_{Y_1}<T_{Y_3} ) P_{\omega}^{Y_1}( T_y<\infty )
      + P_{\omega}^y( T_{Y_3}<T_{Y_1} )P_{\omega}^{Y_3}( T_y<\infty ) )\,,
\end{eqnarray*}
where $G(y,Y_1\land Y_3):=E_{\omega}^y\left[\sum_{k=0}^{T_{Y_1}\land
T_{Y_3}}\ind_{\{X_k=y\}}\right]$. Accordingly,
\begin{eqnarray*}
    E_{\omega}^y[N(y)]={G(y,Y_1\land Y_3)\over
    1- P_{\omega}^y( T_{Y_1}<T_{Y_3} ) P_{\omega}^{Y_1}( T_y<\infty )
      - P_{\omega}^y( T_{Y_3}<T_{Y_1} )P_{\omega}^{Y_3}( T_y<\infty )
      }\,.
\end{eqnarray*}

\noindent Notice that ${\left[
    1- P_{\omega}^y( T_{Y_1}<T_{Y_3} ) P_{\omega}^{Y_1}( T_y<\infty )
      - P_{\omega}^y( T_{Y_3}<T_{Y_1} )P_{\omega}^{Y_3}( T_y<\infty )
      \right]^{-1}}$ is the expected number of times when the walk go
      from $y$ to $Y_1$ or $Y_3$ and then returns to $y$, which is
      naturally smaller than $E_{\omega}^y[N(Y_1)+N(Y_3)]$. We have
\begin{eqnarray*}
E_{\omega}^y[N(Y_1)]&=&P_{\omega}^y(T_{Y_1}<\infty){\left[
1-P_{\omega}^{Y_1}(T_{Y_1}^*<\infty)\right]^{-1}}\\
&\le& {\left[p_1(Y_1,Y_2)\right]^{-1}}\,,
\end{eqnarray*}

\noindent where as before $p_1(Y_1,Y_2) = P_{\omega}^{Y_1}\left(
    T_{{\buildrel \leftarrow \over Y_1}}=\infty\,,\,T_{Y_1}^*=\infty\,,\,T_{Y_2}=\infty
    \right)$. Similarly $E_{\omega}^y[N(Y_3)]\le {\left[
\beta(Y_3)\right]^{-1}}$.  We obtain
\begin{eqnarray}
\label{A1}
    P_{\omega}^e(T_y<\infty)\left(E_{\omega}^y[N(y)]\right)^{\lambda}\le
    {\left[p_1(Y_1,Y_2)\beta(Y_3)\right]^{-1}}P_{\omega}^e(T_y<\infty)\left(G(y,Y_1\land
  Y_3)\right)^{\lambda}\,.
\end{eqnarray}

\noindent We deduce from the Markov property that
$P_{\omega}^e(T_y<\infty)=
P_{\omega}^e(T_{Y_1}<\infty)P_{\omega}^{Y_1}(T_y<\infty) $ and
$P_{\omega}^{Y_1}(T_y<\infty)=G(Y_1,y)P_{\omega}^{Y_1}(T_y<T_{Y_1}^*)$
where
$G(Y_1,y):=E_{\omega}^{Y_1}\left[\sum_{k=0}^{T_y}\ind_{\{X_k=Y_1\}}\right]$.
By Lemma \ref{one}, we have $P_{\omega}^{Y_1}(T_y<T_{{\buildrel
\leftarrow \over Y_1}})\le
\widetilde{P}_{\omega}^{Y_1}(T_y<T_{{\buildrel \leftarrow \over
Y_1}})$. In words, it means that the probability to escape by $y$ is
lower for the RWRE on the tree than for the restriction of the walk
on $[\![Y_1,y]\!]$. Furthermore $G(Y_1,y)\le
E_{\omega}^{Y_1}[N(Y_1)]\le \left[p_1(Y_1,Y_2)\right]^{-1} $, so
that
\begin{eqnarray}
\nonumber
  P_{\omega}^e(T_y<\infty)
&\le&
  {P_{\omega}^e(T_{Y_1}<\infty)\widetilde{P}_{\omega}^{Y_1}(T_y<T_{{\buildrel \leftarrow \over Y_1}})}{\left[p_1(Y_1,Y_2)\right]^{-1}}\\
&\le&
  {P_{\omega}^e(T_{Y_1}<\infty)\left(\widetilde{P}_{\omega}^{Y_1}(T_y<T_{{\buildrel \leftarrow \over Y_1}})\right)^{\lambda}}{\left[p_1(Y_1,Y_2)\right]^{-1}}\label{A2}\,.
\end{eqnarray}

\noindent We observe that
\begin{eqnarray}
\label{g}
   G(y,Y_1\land Y_3)
=  {\left[1-P_{\omega}^{y}(T_y^*<T_{Y_1}\land
T_{Y_3})\right]^{-1}}\,.
\end{eqnarray}

\noindent Call $y_3$ the unique child of $y$ such that $y_3\le Y_3$.
Consequently,
\begin{eqnarray*}
&&    P_{\omega}^{y}(T_{y}^*<T_{Y_1}\land T_{Y_3})\\
&\le& [ 1 - \omega( y,y_3 ) - \omega( y, {\buildrel \leftarrow \over
    y} ) ] + \omega( y , {\buildrel \leftarrow \over y} )
    P_{\omega}^{{\buildrel
    \leftarrow \over y}}( T_{y}<T_{Y_1} ) +
    \omega(y,y_3)P_{\omega}^{y_3}( T_{y}<T_{Y_3} )\,.
\end{eqnarray*}

\noindent By Lemma \ref{one}, we have
\begin{eqnarray*}
  P_{\omega}^{{\buildrel \leftarrow \over
  y}}(T_{y}<T_{Y_1})
&\le&
  \widetilde {P}_{\omega}^{{\buildrel
  \leftarrow \over y}}(T_{y}<T_{Y_1})\,,\\
  P_{\omega}^{y_3}(T_{y}<T_{Y_3}) &\le&
  \widetilde{P}_{\omega}^{y_3}(T_{y}<T_{Y_3})\,.
\end{eqnarray*}

\noindent Equation (\ref{g}) becomes $G(y,Y_1\land Y_3)\le
(\omega(y,y_3)+\omega(y,{\buildrel \leftarrow \over
y}))^{\!\!-1}\widetilde {G}(y,Y_1\land Y_3)$ where $\widetilde
{G}(y,Y_1\land Y_3)$ stands for the expectation of the number of
times the one-dimensional RWRE associated to the pair $(Y_1,Y_3)$ by
Lemma \ref{one} crosses $y$ before reaching $Y_1$ or $Y_3$ when
started from $y$. Since $\nu(y)\le b_0$, there exists a constant
$c_{28}$ such that $(\omega(y,{\buildrel \leftarrow \over y}) +
\omega(y,y_3))^{\!\!-1}\le c_{28}$. It yields
\begin{eqnarray}
\label{A3}
  G(y,Y_1\land Y_3)
\le
  c_{28}\,\widetilde {G}(y,Y_1\land Y_3)\,.
\end{eqnarray}

\noindent Finally, using  (\ref{A2}), (\ref{A3}), and the following
inequality,
\begin{eqnarray*}
\label{part2}
  \widetilde{P}_{\omega}^{Y_1}(T_y<T_{{\buildrel \leftarrow \over Y_1}})\,\widetilde {G}(y,Y_1\land Y_3)
\le
  \widetilde{E}^{Y_1}_{\omega}[T_{{\buildrel \leftarrow \over Y_1}}\land T_{Y_3}]\,,
\end{eqnarray*}

\noindent we get
\begin{eqnarray*}
 P_{\omega}^e(T_y<\infty)\left(G(y,Y_1\land
  Y_3)\right)^{\lambda}\le {c_{28}\over p_1(Y_1,Y_2)}P_{\omega}^e(T_{Y_1}<\infty)(\widetilde{E}_{\omega}^{Y_1}[T_{{\buildrel \leftarrow \over Y_1}}\land
  T_{Y_3}])^{\lambda}\,.
\end{eqnarray*}

\noindent By Lemma \ref{estimation}, for any $y\in \t$, we have
$$
\widetilde{E}_{\omega}^{Y_1}[T_{{\buildrel \leftarrow \over
Y_1}}\land
  T_{Y_3}]\le c_{19}(n_0)\widetilde{E}_{\omega}^{Y_2}[T_{{\buildrel \leftarrow \over Y_2}}\land
  T_{Y_3}]\,.
$$

\noindent It follows that
\begin{eqnarray}
\label{i2} P_{\omega}^e(T_y<\infty)\left(G(y,Y_1\land
  Y_3)\right)^{\lambda}\le {c_{28}c_{19}^{\lambda}\over
  p_1(Y_1,Y_2)}P_{\omega}^e(T_{Y_1}<\infty)(\widetilde{E}^{Y_2}_{\omega}[T_{{\buildrel \leftarrow \over Y_2}}\land
  T_{Y_3}])^{\lambda}\,.
  \end{eqnarray}

 \noindent In view of equations (\ref{A1}) and
 (\ref{i2}), we obtain
\begin{eqnarray*}
  P_{\omega}^e(T_y<\infty)\left(E_{\omega}^y[N(y)]\right)^{\lambda}
&\le&
  c_{29}\,P_{\omega}^e(T_{Y_1}<\infty)H(Y_1,y,Y_3)
\end{eqnarray*}

\noindent where
\begin{eqnarray*}
  H(Y_1,y,Y_3)
:=
  {\left[p_1(Y_1,Y_2)^2\beta(Y_3)\right]^{-1}}
\left(\widetilde{E}^{Y_2}_{\omega}[T_{{\buildrel \leftarrow \over
Y_2}}\land
  T_{Y_3}]\right)^{\lambda}\,.
\end{eqnarray*}

\noindent By equation (\ref{sumgen}), it implies that
$$
E_{\q}[N_{2,n}^{\lambda}] \le
  c_{29}\,E_{\Q}\left[\sum_{|y|= n}
  P_{\omega}^e(T_{Y_1}<\infty)
  H(Y_1,y,Y_3)\right]\,.
$$

\noindent Arguing over the value of $Y_1$ gives
\begin{eqnarray*}
E_{\q}[N_{2,n}^{\lambda}] &\le& c_{29}\,E_{\Q}\left[\sum_{|z|\le
  n-n_0}P_{\omega}^e(T_{z}<\infty)\left(
  \sum_{|y|=n,Y_1=z}H(z,y,Y_3)\right)\right]\\
&=&c_{29}\,E_{\Q}\left[\sum_{|z|\le
  n-n_0}P_{\omega}^e(T_{z}<\infty)E_{\Q}\left[\sum_{|y|=n-|z|,Y_1=e}H(e,y,Y_3)\right]\right]\\
&=&c_{29}\,E_{\Q}\left[\sum_{|z|\le
  n-n_0}P_{\omega}^e(T_{z}<\infty)S(n-|z|)\right]\,,
\end{eqnarray*}

\noindent by equation (\ref{defs}). Lemma \ref{rn} yields that
\begin{eqnarray*}
E_{\q}[N_{2,n}^{\lambda}] &\le&
c_1c_{29}\sum_{k = n_0}^{n}S(k)\\
&\le& c_1c_{29}\sum_{k\ge n_0}S(k)\,.\;\Box
\end{eqnarray*}

We call as before
$m(n,\lambda):=\E\left[\left(E_{\omega}^0\left[T_{-1}\land
T_n\right]\right)^{\lambda}\right]$ for the one-dimensional RWRE
$(R_n)_{n\ge 0}$. The following lemma gives an estimate of $S(n)$.
\begin{lemma}
\label{S} There exists a constant $c_{30}$ such that for any
$\ell\ge 0$,
\begin{eqnarray*}
\label{s1}S(\ell+n_0)\le c_{30}\sum_{i\ge \ell}m(i,\lambda)r^{i}\,.
\end{eqnarray*}
\end{lemma}

\noindent {\it Proof.} Let $\ell\ge 0$ and
$f(Y_2,Y_3):=\left(\widetilde{E}^{Y_2}[T_{{\buildrel \leftarrow\over
Y_2}}\land T_{Y_3}]\right)^{\lambda}$. We have
\begin{eqnarray*}
  S(\ell+n_0)&=&E_{\Q}\left[\sum_{|y|=\ell+n_0:Y_1=e}{\left[p_1(e,Y_2)^2\beta(Y_3)\right]^{-1}}f(Y_2,Y_3)\right]\\
  &=&E_{\Q}\left[\sum_{|u|= n_0}{\left[p_1(e,u)\right]^{-2}}\sum_{|y|=\ell+n_0:Y_2=u} f(u,Y_3)\left[\beta(Y_3)\right]^{-1}\right]\,.
\end{eqnarray*}

\noindent If we call $\t_u$ the subtree of $\t$ rooted in $u$, we
observe that
\begin{eqnarray*}
\sum_{|y|=\ell+n_0:Y_2=u} f(u,Y_3)\left[\beta(Y_3)\right]^{-1} \le
\ind_{\{Z_{n_0}>b_0\}} \sum_{|z|\ge\ell+n_0:z\in\mathbb{W}(\t_u)}
f(u,z)\left[\beta(z)\right]^{-1}\ind_{\{\nu(z,n_0)>b_0\}}\,,
\end{eqnarray*}

\noindent where $\mathbb{W}$ was defined in equation (\ref{defw}).
The Galton--Watson property yields that
\begin{eqnarray*}
\label{B} S(\ell+n_0)&\le& E_{\Q}\left[\sum_{|u|=
n_0}\frac{\ind_{\{Z_{n_0}>b_0\}}}{p_1(e,u)^2}\right]
E_{\Q}\left[\sum_{|z|\ge
\ell,z\in\mathbb{W}}f(e,z)\left[\beta(z)\right]^{-1}\ind_{\{\nu(z,n_0)>b_0\}}\right]\\
&=&E_{\Q}\left[\sum_{|u|=
n_0}\frac{\ind_{\{Z_{n_0}>b_0\}}}{p_1(e,u)^2}\right]E_{\Q}\left[\sum_{|z|\ge
\ell,z\in\mathbb{W}}f(e,z)\right]E_{\Q}\left[{\ind_{\{Z_{n_0}>b_0\}}\over \beta(e)}\right]\\
&\le&c_{31}E_{\Q}\left[\sum_{|z|\ge
\ell,z\in\mathbb{W}}f(e,z)\right]\,,
\end{eqnarray*}

\noindent by Lemma \ref{p1} and equation (\ref{fin}). The proof
follows then from
\begin{eqnarray*}
E_{\Q}\left[\sum_{|z|\ge \ell,z\in\mathbb{W}}f(e,z)\right] &=&
E_{GW}\left[\sum_{|z|\ge
\ell,z\in\mathbb{W}}m(|z|,\lambda)\right]\\
&=& \sum_{i:in_0\ge \ell}m(in_0,\lambda)E_{GW}[W_i]\le \sum_{in_0\ge
\ell}m(in_0,\lambda)r^{in_0} \,,
\end{eqnarray*}

\noindent where the last inequality comes from equation (\ref{W}).
$\Box$\\

We are now able to prove (\ref{N2}).\\
{\it Proof of Lemma \ref{lem}, equation (\ref{N2})}. By Lemma
\ref{N31},
$$
E_{\q}[N_{2,n}^{\lambda}] \le c_{27}\sum_{\ell\ge 0}S(\ell+n_0)\,.
$$
\noindent Lemma \ref{S} tells that
\begin{eqnarray*}
\sum_{\ell\ge 0}S(\ell+n_0) &\le&
  c_{30}\sum_{i\ge \ell \ge 0}m(i,\lambda)r^i\\
&=&
  c_{30}\sum_{i\ge 0}(i+1)m(i,\lambda)r^i\,,
\end{eqnarray*}

\noindent which is finite by equation (\ref{lambda}). $\Box$

\section{Proof of Theorem \ref{global}}

If we suppose that $\Lambda<1$, then Theorem \ref{kesten} ensures
that ${|X_n|\over n}$ tends to $0$. Suppose now that $\Lambda>1$.
Take $\lambda=1$ in the proof of the lower bound of Theorem
\ref{kesten} in Section 5 to see that $|X_n|\ge {n\over c_{23}}$ for
sufficiently large n, which proves the positivity of the speed in
this case. Theorem \ref{global} is proved. $\Box$

\section{Proof of Theorem \ref{reinforce}}

When $b\ge 3$, Theorem \ref{reinforce} follows immediately from
Theorem \ref{errw}. In the rest of this section, we assume that $\t$
is a binary tree. Thanks to the correspondence between RWRE and
LERRW mentioned in the introduction, we only have to prove the
positivity of the speed for a RWRE on the binary tree such that the
density of $\omega(y,z)$ on $(0,1)$ is given by
\begin{eqnarray}
\label{om1} f_0(x) &=& 1 ~~\qquad \qquad \qquad \qquad \qquad \mbox{if}~z={\buildrel \leftarrow \over y}\,\\
\label{om2} f_1(x) &=& {1\over \Gamma({1\over 2})\Gamma({3\over
2})}x^{-1/2}(1-x)^{1/2}~~\mbox{if}~z~\mbox{is a child of}~y.
\end{eqnarray}

We propose to prove three lemmas before handling the proof of the
theorem.
\begin{lemma}
\label{beta3} We have for any $0<\delta<1$, $$\E\left[{1\over
\beta^{\delta}}\right]<\infty\,.
$$
\end{lemma}
{\it Proof.} By equation (\ref{beta2}), for any $y\in \t$,
\begin{eqnarray*}
{1 \over \beta(y)^{\delta}} &\le& \left(1 + \min_{i=1,2}{1\over
A(y_i)\beta(y_i)} \right)^{\delta}\\
&\le& 1 + \min_{i=1,2}{1 \over
A(y_i)^{\delta}\beta(y_i)^{\delta}}\,.
\end{eqnarray*}

\noindent Notice that by (\ref{om1}),
$$
\E\left[\min_{i=1,2}{1\over A(y_i)^{\delta}}\right]\le
2^{\delta}\,\E\left[\left({1\over
A(y_1)+A(y_2)}\right)^{\delta}\right]=2^{\delta}\,\E\left[\left({\omega(y,\buildrel
\leftarrow \over y)\over 1-\omega(y,\buildrel \leftarrow \over
y)}\right)^{\delta}\right]<\infty\,.
$$

\noindent The proof is therefore the proof of Lemma \ref{beta} when
replacing $A(y)$ and $\beta(y)$ respectively by $A(y)^{\delta}$ and
$\beta(y)^{\delta}$. $\Box$ \\

Recall that for any $y\in\t$, $\gamma(y):=P_{\omega}^y(T_{{\buildrel
\leftarrow \over y}}=\infty,\,T_y^*=\infty)$.
\begin{lemma}
\label{reinfor2} There exists $\mu \in(0,1) $ such that for any
$\varepsilon\in(0,1)$, we have
$$\E\left[\left({\ind_{\{ \omega(e,\buildrel \leftarrow \over e)\le
1-\varepsilon\}}\over \gamma(e)}\right)^{1/\mu}\right]<\infty\,.
$$
\end{lemma}
{\it Proof.} We see that
\begin{eqnarray*}
{1\over \gamma(e)} = {1\over
\omega(e,e_1)\beta(e_1)+\omega(e,e_2)\beta(e_2)} \le
\min_{i=1,2}{1\over \omega(e,e_i)\beta(e_i)}\,.
\end{eqnarray*}

\noindent Let $\mu\in(0,1)$ and $\varepsilon\in(0,1)$. We compute
$\P(\omega(e,\buildrel \leftarrow \over e)\le
1-\varepsilon\,,\,\min_{i=1,2}\,\left\{[{
\omega(e,e_i)\beta(e_i)}]^{-1/\mu}\right\}>n)$ for $n\in \r_+^*$. We
observe that $\{\omega(e,\buildrel \leftarrow \over e)\le
1-\varepsilon\}\subset \{\omega(e,e_1)\ge
\varepsilon/2\}\cup\{\omega(e,e_2)\ge \varepsilon/2\}$. By symmetry,
\begin{eqnarray*}
&&\P\left(\omega(e,\buildrel \leftarrow \over e)\le
1-\varepsilon\,,\,\min_{i=1,2}\,
\left\{[\omega(e,e_i)\beta(e_i)]^{-1/\mu}\right\}>n\right)\\
 &\le&
2\P\left(\omega(e,e_2)\ge \varepsilon/2\,,\, \min_{i=1,2}\,
\left\{[\omega(e,e_i)\beta(e_i)]^{-1/\mu}\right\}>n\right)\\
&\le& 2\P\left(\beta(e_2)^{-1}>
n^{\mu}\varepsilon/2\,,\,[\omega(e,e_1)\beta(e_1)]^{-1/\mu}>n\right)\\
&\le& 2\P\left(\beta(e_2)^{-1}>
n^{\mu}\varepsilon/2\,,\,\omega(e,e_1)\le n^{-1/2}\right) +
2\P\left(\beta(e_2)^{-1}> n^{\mu}\varepsilon/2\,,\,\beta(e_1)^{-1}>
n^{\mu-1/2}\right)\\
&=:& 2\P(E_1)+2\P(E_2)\,.
\end{eqnarray*}

\noindent Let $0<\delta<1$. We have by (\ref{om2}) and Lemma
\ref{beta3},
\begin{eqnarray*}
\P(E_1)&=&\P\left(\omega(e,e_1)\le n^{-1/2})\P({\beta(e_2)^{-1}}>
n^{\mu}\varepsilon/2\right)\\
&\le& c_{32}n^{-1/4}n^{-\delta\mu}\,.
\end{eqnarray*}

\noindent Similarly,
\begin{eqnarray*}
\P(E_2)  &=& \P\left({\beta(e_1)^{-1}}> n^{\mu-1/2})\P({
\beta(e_2)^{-1}}>
n^{\mu}\varepsilon/2\right)\\
&\le& c_{33}n^{-\delta(\mu-1/2)}n^{-\delta\mu}\,.
\end{eqnarray*}

\noindent It suffices to take $1/4+\delta\mu>1$ and
$\delta(2\mu-1/2)>1$ to complete the proof, for example by taking
$\delta=4/5$ and $\mu=19/20$. $\Box$ \\

Let $\varepsilon \in (0,1/3)$ be such that
\begin{equation}
\label{var} \E\left[\left(\#\{i\,:\,\omega(e_i,e)>
1-\varepsilon\}\right)^{{2-\mu\over 1-\mu}}\right]<1\,.
\end{equation}

\noindent Denote by $\mathbb{U}$ the set of the root and all the
vertices $y$ such that for any vertex $x\in \t$ with $e<x\le y$, we
have $\omega(x,\buildrel \leftarrow \over x)> 1-\varepsilon$; we
observe that by (\ref{var}), $\mathbb{U}$ is a subcritical
Galton--Watson tree. Denote by $U_k$ the size of the generation $k$.

\begin{lemma}
\label{w2} There exists a constant $c_{34}<1$ such that for any
$k\ge 0$
$$
\E\left[U_k^{1/(1-\mu)}\right] \le c_{34}^k\,.
$$
\end{lemma}
{\it Proof.} By Galton--Watson property,
\begin{eqnarray*}
\E\left[U_{k+1}^{1/(1-\mu)}\right]=\E\left[\left(\sum_{i=1}^{U_1}U_k^{(i)}\right)^{1/(1-\mu)}\right]
\end{eqnarray*}

\noindent where conditionally on $U_1$, $U_k^{(i)},\,i\ge 1$ is a
family of i.i.d random variables distributed as $U_k$. Since
$\left(\sum_{i=1}^n a_i\right)^p \le n^p\sum_{i=1}^pa_i^p$ (for
$p>0$ and $a_i\ge 0$), it yields that
\begin{eqnarray*}
\E\left[U_{k+1}^{1/(1-\mu)}\right] &\le& \E\left[U_1^{1/(1-\mu)}\sum_{i=1}^{U_1}\left(U_k^{(i)}\right)^{1/(1-\mu)}\right]\\
&=& \E\left[U_1^{{2-\mu\over
1-\mu}}\right]\E\left[U_k^{1/(1-\mu)}\right]\,.
\end{eqnarray*}

\noindent The proof follows from equation (\ref{var}). $\Box$ \\

We are now able to complete the proof of Theorem \ref{reinforce}.\\

\noindent{\it Proof of Theorem \ref{reinforce} : the binary tree
case.} We suppose without loss of generality that
$\omega(e,{\buildrel \leftarrow \over e})\le 1-\varepsilon$. For any
vertex $y$, we call $Y$ the youngest ancestor of $y$ such that
$\omega(Y,\buildrel \leftarrow \over Y)\le 1-\varepsilon$. We have
for any $n\ge 0$,
\begin{eqnarray*}
E_{\omega}^e[N_n] =
\sum_{|y|=n}P_{\omega}^e(T_y<\infty)E_{\omega}^y[N(y)]\,,
\end{eqnarray*}

\noindent where, as before, $N(y):=\sum_{k\ge 0}\ind_{\{X_k=y\}}$
and $N_n=\sum_{|y|=n}N(y)$. By the Markov property,
\begin{eqnarray*}
E_{\omega}^y[N(y)] &=& G(y,Y) +
P_{\omega}^y(T_{Y}<\infty)P_{\omega}^{Y}(T_y<\infty)E_{\omega}^{y}[N(y)]\,,
\end{eqnarray*}

\noindent where
$G(y,Y):=E_{\omega}^y\left[\sum_{k=0}^{T_{Y}}\ind_{\{X_k=y\}}\right]$.
It yields that
\begin{eqnarray*}
E_{\omega}^e[N_n]&=&
\sum_{|y|=n}P_{\omega}^e(T_y<\infty){G(y,Y)\over
1-P_{\omega}^{Y}(T_y<\infty)P_{\omega}^y(T_{Y}<\infty)}\\
&\le& \sum_{|y|=n}P_{\omega}^e(T_y<\infty){G(y,Y)\over
1-P_{\omega}^{Y}(T_{Y}^*<\infty)}\\
&\le& \sum_{|y|=n}P_{\omega}^e(T_y<\infty){G(y,Y)\over \gamma(Y)}\,.
\end{eqnarray*}

\noindent By coupling the walk on $[\![y,Y]\!]$ with a
one-dimensional random walk, we see that
$P_{\omega}^{y}(T^*_y<T_{Y})\le \varepsilon +
(1-\varepsilon){\varepsilon\over 1-\varepsilon}=2\varepsilon \le
2/3$, so that $G(y,Y)\le 3$. On the other hand,
$P_{\omega}^e(T_y<\infty)\le P_{\omega}^e(T_{Y}<\infty)$. Therefore,
\begin{eqnarray*}
\e[N_n] &\le& 3\E\left[\sum_{|y|=n}P_{\omega}^e(T_{Y}<\infty){1\over
\gamma(Y)}\right]\\
&=&3\E\left[\sum_{|y|=n}\sum_{z=Y}P_{\omega}^e(T_{z}<\infty){1\over
\gamma(z)}\right]\\
&=& 3 \E\left[\sum_{|z|\le
n}P_{\omega}^e(T_{z}<\infty)\sum_{|y|=n:Y=z}{1\over
\gamma(z)}\right]\,.
\end{eqnarray*}

\noindent By independence and stationarity of the environment,
\begin{eqnarray*}
\e[N_n]&\le& 3\sum_{|z|\le n}
\p(T_z<\infty)\E\left[\sum_{|y|=n-|z|:Y=e}{1\over \gamma(e)}\right]\\
&=& 3\sum_{|z|\le n}
\p(T_z<\infty)\E\left[{\ind_{\{\omega(e,{\buildrel \leftarrow \over
e})\le 1-\varepsilon\}}U_{n-|z|}\over \gamma(e)}\right]\\
&\le& 3\sum_{|z|\le
n}\p(T_z<\infty)\E\left[\left({\ind_{\omega(e,\buildrel\leftarrow\over
e)\le 1-\varepsilon}\over
\gamma(e)}\right)^{1/\mu}\right]^{\mu}\E\left[U_{n-|z|}^{1/(1-\mu)}\right]^{1-\mu}\,,
\end{eqnarray*}

\noindent by the H\"older inequality. We use Lemmas \ref{reinfor2}
and \ref{w2} to see that
\begin{eqnarray*}
\e[N_n] \le c_{35}\sum_{|z|\le n}\p(T(z)<\infty)c_{36}^{n-|z|}.
\end{eqnarray*}

\noindent By Lemma \ref{rn},
\begin{eqnarray*}
\e[N_n] \le c_{35}c_{1}\sum_{k=0}^nc_{36}^k <
c_{35}c_1/(1-c_{36})\,.
\end{eqnarray*}

\noindent Since $\tau_n \le \sum_{k=-1}^nN_k$ and $N_{-1}\le N_0$,
where $\tau_n:=\inf\left\{k\ge 0\, : |X_k|=n\right\}$ as before, we
have $\e[\tau_n]\le c_{37}\,n$. Fatou's lemma yields that
$\p$-almost surely, $ \liminf_{n\rightarrow \infty}{\tau_n \over
n}<\infty $, which proves that $v>0$ in view of the relation
$\lim_{n\rightarrow\infty}{\tau_n\over n}={1\over v}$. $\Box$

\section{Proof of Lemmas \ref{tec1} and \ref{infeq2}}

We consider the one-dimensional RWRE $(R_n)_{n\ge 0}$ when we
consider the case $\t=\{-1,0,1,\ldots\}$. This RWRE is such that the
random variables $A(i)$, $i\ge 0$ are independent and have the
distribution of $A$, when we set for $i\ge 0$,
$$
A(i):={\omega(i,i+1)\over \omega(i,i-1)}
$$
with $\omega(y,z)$ the quenched probability to jump from $y$ to $z$.
We recall that, as defined in equations (\ref{pna}) and (\ref{mn}),
\begin{eqnarray*}
p\,(n,a)&:=&\p^0(T_{-1}\land T_n>a)\,,\\
m(n,\lambda)&:=&\E\left[\left(E_{\omega}^0\left[T_{-1}\land
T_n\right]\right)^{\lambda}\right]\,.
\end{eqnarray*}

\noindent We study the walk $(R_n)_{n\ge 0}$ through its potential.
We introduce for $p\ge i\ge 0$, $V(0)=0$ and
\begin{eqnarray*}
  V(i)    &=& -\sum_{k=0}^{i-1}\ln(A(k))\,,\\
  M(i)    &=& \max_{0\le k\le i} V(k)\,,\\
  H_1(i)  &=& \max_{0\le k\le i} V(k)-V(i)\,,\\
  H_2(i,p)&=& \max_{i\le k\le p} V(k)-V(i)\,.
\end{eqnarray*}

\noindent Let us introduce for $t\in \r$ the Laplace transform
$\E[A^t]$, and define $\phi(t):=\ln(\E[A^t])$. Denote by $I$ its
Legendre transform $I(x)=\sup\{tx-\phi(t),t\in\r\}$ where $x\in\r$.
Let also
$$
[a,b]:=[\mbox{ess inf}(\ln A),\mbox{ess sup}(\ln A)]\,.
$$
 Two situations occur. If $a=b$, it means that $A$ is a constant almost surely. In this case, $I(x)=0$ if $x=a$ and is infinite otherwise.
 If $a<b$, then $I$ is finite on $]a,b[$ and infinite on $\r\backslash[a,b]$. Moreover, for any $x\in]a,b[$, we have $I'(x)=t(x)$
 where $t(x)$ is the real
 such that $I(x)=xt(x)-\phi(t(x))$, or, equivalently, $x=\phi'(t(x))$.

We define and compute two useful parameters. Call
$\mathcal{D}:=\{x_1,\,x_2,\,,z_1,\,z_2\in \r_+^4, z_1+z_2\le 1\}$.
Define for $0<\lambda\le 1$, and with the convention that
$0\times\infty:=0$,
\begin{eqnarray}
\label{L1}   L(\lambda) &:=&
  \sup_{\mathcal{D}}\bigg\{\bigg((x_1z_1)\land (x_2z_2)\bigg)\lambda-I(-x_1)z_1-I(x_2)z_2\bigg\}\,,\\
\label{l'} L' &:=& \sup\bigg\{{x_1+x_2\over
x_1x_2}\ln(q_1)-{I(-x_1)\over x_1
  }-{I(x_2)\over x_2}\,,\, x_1,x_2> 0\bigg\}\,.
\end{eqnarray}

\noindent If $q_1=0$, we set $L'=-\infty$. Notice that
$L(\lambda)\ge 0$ is necessarily reached for $x_1z_1=x_2z_2$. It
yields that
\begin{eqnarray}
\label{L2}
  L(\lambda)=
  0\vee \sup \bigg\{{x_1x_2\over x_1+x_2}\lambda-I(-x_1){x_2\over x_1
  +x_2}-I(x_2){x_1\over x_1+x_2}\,,\, x_1,x_2> 0\bigg\}\,,
\end{eqnarray}

\noindent where $c\lor d:=\max(c,d)$. The computation of
$L(\lambda)$ and $L'$ is done in the following lemma.

\begin{lemma}
\label{est2} We have
\begin{eqnarray}
\label{es21} L(\lambda)&=&0\vee \phi(\bar t\,)\,,\\
\label{es22} L'&=&-\Lambda\,,
\end{eqnarray} where $\bar t$ verifies $\phi(\bar t\,)=\phi(\bar t+\lambda)$ if it exists and $\bar t :=0$ otherwise.
\end{lemma}

\noindent {\it Proof.} When $A$ is a constant almost surely,
$L(\lambda)=0$ and (\ref{es21}) is true. Therefore we assume that
$a<b$. Considering equation (\ref{L2}), we see that if
$L(\lambda)>0$, then $L(\lambda)$ is reached by a pair $(x_1,x_2)$
which
 satisfies:
\begin{eqnarray}
\label{e1}
  \lambda {x_2\over x_1+x_2}+{I(-x_1)\over x_1+x_2}+I'(-x_1)-{I(x_2)\over x_1+x_2}
  &=&
  0\,,\\
\label{e2}
   \lambda {x_1\over x_1+x_2}-{I(-x_1)\over x_1+x_2}+{I(x_2)\over x_1+x_2}-I'(x_2)
&=&
  0\,.
\end{eqnarray}

\noindent We deduce from equations (\ref{e1}) and (\ref{e2}) that
$I'(x_2)-I'(-x_1)=\lambda$, i.e. $t(x_2)-t(-x_1)=\lambda$.  Plugging
this into (\ref{L2}) yields
\begin{eqnarray*}
  L(\lambda)=
  0\vee\sup\left\{{\phi(t)\phi'(t+\lambda)-\phi(t+\lambda)\phi'(t)\over
  \phi'(t+\lambda)-\phi'(t)},t\in\r,\phi'(t)< 0,\phi'(t+\lambda)> 0\right\}\,.
\end{eqnarray*}

\noindent Let
$h(t):={\phi(t)\phi'(t+\lambda)-\phi(t+\lambda)\phi'(t)\over
\phi'(t+\lambda)-\phi'(t)}$. Then $L(\lambda)=0\vee h(\bar t\,)$
where $\bar t$ verifies $h'(\bar t\,)=0$, which is equivalent to say
that $\phi(\bar t\,)=\phi(\bar t +\lambda)$. We find that $h(\bar
t\,)=\phi(\bar t\,)$, which gives (\ref{es21}). The computation of (\ref{es22}) is similar and is therefore omitted. $\Box$\\

\subsection{Proof of Lemma \ref{tec1}}

We begin by some notation. Let $A>0$ and $B>0$ be two expressions
which can depend on any variable, and in particular on $n$. We say
that $ A\lesssim B$ if we can find a function $f$ of the variable
$n$ such that $\lim_{n\rightarrow\infty}{1\over n}\ln(f(n))=0$ and
$A\le f(n)B$. We say that $A\simeq B$ if $ A\lesssim B$ and $
B\lesssim A$. By circuit analogy (see \cite{ds84}), we find for
$0\le i\le n$,
\begin{eqnarray*}
P_{\omega}^0\left(T_i<T_{-1}\right)=\frac{1}{e^{V(0)}+e^{V(1)}+
\ldots + e^{V(i)}}\,.
\end{eqnarray*}

\noindent It follows that
\begin{eqnarray}
\label{cir1} {e^{-M(i)}\over n+1}\le
P_{\omega}^0\left(T_i<T_{-1}\right)\le e^{-M(i)}\,.
\end{eqnarray}

\noindent We deduce also that
\begin{eqnarray}
\label{cir2}
{e^{-H_2(i,n)}\over n+1}\le P_{\omega}^{i+1}\left(T_{n}<T_{i}\right)\le e^{-H_2(i,n)}\,,\\
{e^{-H_1(i)}\over n+1}\le
P_{\omega}^{i-1}\left(T_{-1}<T_{i}\right)\le
e^{-H_1(i)}\,.\label{cir3}
\end{eqnarray}

\noindent Finally, the quenched expectation $G\left(i,-1\land
n\right)$ of the number of times the walk starting from $i$ returns
to $i$ before reaching $-1$ or $n$ verifies
\begin{eqnarray*}
G\left(i,-1\land
n\right)&=&\left\{\omega(i,i-1)P_{\omega}^{i-1}\left(T_{-1}<T_{i}\right)+\omega(i,i+1)P_{\omega}^{i+1}\left(T_{n}<T_{i}\right)\right\}^{-1}\,,
\end{eqnarray*}

\noindent so that
\begin{eqnarray*}
c_{37}e^{H_1(i)\land H_2(i,n)}\le G(i,-1\land n)\le
c_{38}(n+1)e^{H_1(i)\land H_2(i,n)}\,.
\end{eqnarray*}

\noindent Since $E_{\omega}^0[T_{-1}\land
T_n]=1+\sum_{i=0}^{n-1}P_{\omega}^0\left(T_i<T_{-1}\right)G\left(i,-1\land
n\right)$, we get
\begin{eqnarray*}
1+{c_{37}\over n+1}\max_{0\le i\le n}e^{-M(i)+H_1(i)\land
H_2(i,n)}\le E_{\omega}^0[T_{-1}\land T_n]\le
1+c_{38}n(n+1)\max_{0\le i\le n}e^{-M(i)+H_1(i)\land H_2(i,n)}\,.
\end{eqnarray*}

\noindent As a result,
\begin{eqnarray}
\label{est}
  \E[\left(E_{\omega}^0[T_{-1}\land T_n]\right)^{\lambda}]
\simeq
  \max_{0\le i\le n}\E\left[e^{\lambda\left[-M(i)+H_1(i)\land H_2(i,n)\right]}\right]\,.
\end{eqnarray}

\noindent We proceed to the proof of Lemma \ref{tec1}. Let $\eta>0$
and $0\le i\le n$. Let $\varepsilon>0$ be such that $(|a|\lor
|b|)\varepsilon<\eta$. For fixed $i$ and $n$, we denote by $K_1$ and
$K_2$ the integers such that
\begin{eqnarray*}
        K_1\eta \le &H_1(i)& <(K_1+1)\eta \,, \\
        K_2\eta \le &H_2(i,n)& < (K_2+1)\eta \,.
\end{eqnarray*}

\noindent Similarly, let $L_1$ and  $L_2$ be integers such that
\begin{eqnarray*}
\exists ~~L_1\lfloor\varepsilon
        n\rfloor \le x< (L_1+1)\lfloor\varepsilon n\rfloor~~ &\mbox{such that}&~~H_1(i)=V(i-x)-V(i)\,,\\
\exists ~~L_2\lfloor\varepsilon n\rfloor \le y
<(L_2+1)\lfloor\varepsilon
        n\rfloor~~ &\mbox{such that}&~~H_2(i,n)=V(i+y)-V(i)\,.
\end{eqnarray*}

\noindent Finally, $e^{\lambda\left[-M(i)+H_1(i)\land
H_2(i,n)\right]}\le e^{(K_1\land
  K_2+1)\lambda \eta n}$. By our choice of $\varepsilon$, we have for any integers $k_1,k_2,\ell_1,\ell_2$,
\begin{eqnarray*}
  \p\left(K_1=k_1,\,L_1=\ell_1\right)
&\le&
  \p\bigg(V(\ell_1\lfloor \varepsilon
n\rfloor)\in\left[-(k_1+2)\eta n,-(k_1-1)\eta n\right]\bigg)\,,
\\
  \p\left(K_2=k_2,\,L_2=\ell_2\right)
&\le&
  \p\bigg(V\left(\ell_2\lfloor \varepsilon n\rfloor\right)\in\left[(k_2-1)\eta n,(k_2+2)\eta n\right]\bigg)\,.
\end{eqnarray*}

\noindent By Cram\'er's theorem (see \cite{holbook} for example),
\begin{eqnarray*}
  \p\bigg(V(\ell_1\lfloor \varepsilon n\rfloor)\in\left[-(k_1+2)\eta n,-(k_1-1)\eta n\right]\bigg)
&\lesssim&
  \exp\bigg(-\ell_1\lfloor \varepsilon n\rfloor(I(-x_1)-\lambda\eta)\bigg)
\\
  \p\bigg(V(\ell_2\lfloor \varepsilon n\rfloor)\in\left[(k_2-1)\eta n,(k_2+2)\eta n\right]\bigg)
&\lesssim&
  \exp\bigg(-\ell_2\lfloor \varepsilon n\rfloor(I(x_2)-\lambda\eta)\bigg)
\end{eqnarray*}

\noindent if $-x_1$ is the point of  $\bigg[{-(k_1+2)\eta n\over
\ell_1\lfloor \varepsilon n\rfloor},{-(k_1-1)\eta n\over
\ell_1\lfloor \varepsilon n\rfloor}\bigg]$ where $I$ reaches the
minimum on this interval, and $x_2$ is the equivalent in
$\bigg[{(k_2-1)\eta n\over \ell_2\lfloor \varepsilon
n\rfloor},{(k_2+2)\eta n\over \ell_2\lfloor \varepsilon
n\rfloor}\bigg]$. It yields that
\begin{eqnarray*}
  &&\e\left[e^{\lambda\left[-M(i)+H_1(i)\land H_2(i,n)\right]}\right]\\
&\lesssim&
 \max_{k_1,k_2,\ell_1,\ell_2\in D'} \exp\left(\left(k_1\land k_2\right)\lambda \eta n -I(-x_1)\ell_1\lfloor \varepsilon n\rfloor-I(x_2)\ell_2\lfloor \varepsilon n\rfloor+3\lambda\eta n\right)\,,
\end{eqnarray*}

\noindent where $D'$ is the (finite) set of all possible values of
$(K_1,K_2,L_1,L_2)$. We note that
\begin{eqnarray*}
  &&(k_1\land k_2)\lambda \eta n -I(-x_1)\ell_1\lfloor \varepsilon
  n\rfloor-I(x_2)\ell_2\lfloor \varepsilon n\rfloor\\
&\le&
  (x_1\ell_1\lfloor \varepsilon n\rfloor \land x_2\ell_2\lfloor \varepsilon
  n\rfloor)\lambda-I(-x_1)\ell_1\lfloor \varepsilon n\rfloor-I(x_2)\ell_2\lfloor \varepsilon n\rfloor+3\lambda\eta n\\
&\le&
  (L(\lambda)+3\lambda\eta)n
\end{eqnarray*}

\noindent by (\ref{L1}). Finally, $\e[e^{\lambda(-M(i)+H_1(i)\land
H_2(i,n))}] \lesssim e^{n(L(\lambda)+6\lambda\eta)}$ so that, by
equation (\ref{est}), $m(n,\lambda)\lesssim
e^{n(L(\lambda)+6\lambda\eta)}$. We let $\eta$ tend to $0$ to get
that
$$
\limsup_{n\rightarrow\infty} {1\over n}\ln(m(n,\lambda))\le
L(\lambda)\,.
$$
Let $\lambda<\Lambda$. By definition of $\Lambda$ and equation
(\ref{es21}), it implies that $L(\lambda)<{1\over q_1}$, so that we
can find $r>q_1$ such that $\sum_{n\ge 0}m(n,\lambda)r^n<\infty$. It
means that $\lambda\le \lambda_c$. Consequently, $\Lambda\le
\lambda_c$. $\Box$

\subsection{Proof of Lemma \ref{infeq2}}

Fix $x_1,\,x_2> 0$. Write
\begin{eqnarray*}
  z_1=
  {x_2\over x_1+x_2},~~z_2={x_1\over x_1+x_2},~~z={x_1x_2\over x_1+x_2}\,.
\end{eqnarray*}

\noindent Let $a\ge 100$ and $n=n(a):=\lfloor {\ln(a)\over z}
\rfloor$. We have, by the strong Markov property,
$P_{\omega}^0(T_{-1}\land T_n>a)\ge P_{\omega}^0(T_{\lfloor
z_1n\rfloor}<T_{-1})P_{\omega}^{\lfloor z_1n\rfloor} (T_{\lfloor
z_1n\rfloor}<T_{-1}\land T_n)^a$. It follows by (\ref{cir1}),
(\ref{cir2}) and (\ref{cir3}) that
\begin{eqnarray*}
  p\,(n,a)
&\gtrsim&
  \E\bigg[e^{-M(\lfloor z_1n\rfloor)}\bigg(1-e^{-H_1(\lfloor z_1n\rfloor)\land
  H_2(\lfloor z_1n\rfloor,n)}\bigg)^a\bigg]\\
&\ge&
  (1-e^{-zn})^a\P\bigg(V\left({\lfloor
  z_1n\rfloor}\right)<-zn,\,M\left(\lfloor z_1n\rfloor\right)\le
  0\bigg)\P\bigg(V\left({\lfloor z_2n\rfloor}+1\right)>zn\bigg)\\
&\gtrsim&
  \P\bigg(V\left({\lfloor
  z_1n\rfloor}\right)<-zn,\,M\left(\lfloor z_1n\rfloor\right)\le
  0\bigg)\P\bigg(V\left({\lfloor z_2n\rfloor}+1\right)>zn\bigg)
\end{eqnarray*}

\noindent by our choice of $n$. Let $k\ge 0$. Call $\tau$ the first
time when the walk $(V(i))_{i\ge 0}$ reaches its maximum on $[0,k]$.
Let $i\in [0,k]$ and for $0\le r\le k-1$, $X_r:=\ln(A_{\bar{r}})$
where $\bar{r}:= i+r$ modulo $k$. We observe that
\begin{eqnarray*}
  \P(V_k<-zn,\,\tau=i)
&\le&
  \P(X_0+\ldots+X_{k-1}<-zn,\,X_0+\ldots+X_j\le 0~~ \forall~~ 0\le j\le k-1)\\
&=&
  \P(V_k<-zn,\,M_k\le 0)\,.
\end{eqnarray*}

\noindent We obtain that $\P\left(V_k<-zn,M_k\le 0\right)\ge {1\over
k+1}\P\left(V_k<-zn\right) $. Therefore, for any $\varepsilon>0$,
\begin{eqnarray*}
  p\,(n,a)
&\gtrsim&
  \P\bigg(V\left({\lfloor z_1n\rfloor}\right)<-zn\bigg)\P\bigg(V\left({\lfloor z_2n\rfloor}+1\right)>zn\bigg)\\
&\gtrsim&
  \exp\bigg(n\left(-I(-x_1)z_1-I(x_2)z_2-2\varepsilon\right)\bigg)
\end{eqnarray*}

\noindent by Cram\'er's theorem. It yields that
\begin{eqnarray*}
\liminf_{a\rightarrow\infty}\bigg\{\sup_{\ell\ge 0}
{\ln(q_1^{\ell}p\,(\ell,a))\over \ln(a)}\bigg\}&\ge& \liminf_{a\rightarrow\infty}{\ln(q_1^np\,(n,a))\over \ln(a)}\\
&\ge& {\ln(q_1)-I(-x_1)z_1-I(x_2)z_2-2\varepsilon \over z}\,.
\end{eqnarray*}

\noindent Finally, by (\ref{l'}) and (\ref{es22}),
\begin{eqnarray*}
\liminf_{a\rightarrow\infty}\left\{\sup_{n\ge 0}
{\ln(q_1^np(n,a))\over \ln(a)}\right\} \ge
  L'=-\Lambda\,.\;\Box
\end{eqnarray*}

\noindent {\bf Acknowledgements}: I would like to thank Zhan Shi for
suggesting me the problem and for many precious discussions. I would
also like to thank the referees for their helpful comments.

\bibliographystyle{plain}
\bibliography{biblio}

\bigskip
\bigskip


{\footnotesize

\baselineskip=12pt \hskip60pt  Elie Aid\'ekon

\hskip60pt Laboratoire de Probabilit\'es et Mod\`eles Al\'eatoires

\hskip60pt Universit\'e Paris VI

\hskip60pt 4 Place Jussieu

\hskip60pt F-75252 Paris Cedex 05

\hskip60pt France

\hskip60pt {\tt elie.aidekon@ccr.jussieu.fr}

}

\end{document}